# CLUSTERING IN A STOCHASTIC MODEL OF ONE-DIMENSIONAL GAS

By Vladislav V. Vysotsky[1]

*St. Petersburg State University*

We give a quantitative analysis of clustering in a stochastic model of one-dimensional gas. At time zero, the gas consists of $n$ identical particles that are randomly distributed on the real line and have zero initial speeds. Particles begin to move under the forces of mutual attraction. When particles collide, they *stick* together forming a new particle, called *cluster*, whose mass and speed are defined by the laws of conservation.

We are interested in the asymptotic behavior of $K_n(t)$ as $n \to \infty$, where $K_n(t)$ denotes the number of clusters at time $t$ in the system with $n$ initial particles. Our main result is a functional limit theorem for $K_n(t)$. Its proof is based on the discovered *localization property* of the aggregation process, which states that the behavior of each particle is essentially defined by the motion of neighbor particles.

## 1. Introduction.

1.1. *Description of the model.* We give a quantitative analysis of clustering in a stochastic model of one-dimensional gas. At time zero, the gas consists of $n$ point particles, each one of mass $\frac{1}{n}$. These particles are randomly distributed on the real line and have zero initial speeds. Particles begin to move under the forces of mutual attraction. When two or more particles collide, they *stick* together forming a new particle, called *cluster*, whose mass and speed are defined by the laws of mass and momentum conservation. Between collisions, particles move according to the laws of Newtonian mechanics.

We suppose that the force of mutual attraction does not depend on distance and equals the product of masses. This assumption is natural for

Received March 2007; revised September 2007.
[1]Supported in part by the Grants NSh-4222.2006.1 and DFG-RFBR 436 RUS 113/773/0-1(R).
*AMS 2000 subject classifications.* Primary 60K35, 82C22; secondary 60F17, 70F99.
*Key words and phrases.* Sticky particles, particle systems, gravitating particles, number of clusters, aggregation, adhesion.

one-dimensional models because, by the Gauss law applied to flux of the gravitational field, gravitation is proportional to the distance to the power one minus dimension of the space. At any moment, the acceleration of a particle is thus equal to difference of masses located to the right and to the left of the particle.

Random initial positions of particles are usually described (see [8, 16, 25]) by the following natural models: in the *uniform* model, $n$ particles are independently and uniformly spread on $[0, 1]$; in the *Poisson* model, particles are located at points $\frac{1}{n}S_1, \frac{1}{n}S_2, \ldots, \frac{1}{n}S_n$, where $S_i$ is a standard exponential random walk. In other words, particles are located at points of first $n$ jumps of a Poisson process with intensity $n$.

These two models are the most natural and interesting; let us call them the *main models* of initial positions. However, we will see that behavior of the Poisson model is essentially defined by independence of initial distances between particles rather than by the particular type of the distances' distribution. Therefore, it is of a great mathematical interest to generalize the Poisson model by introducing the *i.d. model*, where "i.d." stands for "independent distances," as follows. Particles are initially located at $\frac{1}{n}S_1, \frac{1}{n}S_2, \ldots, \frac{1}{n}S_n$, where $S_i$ is a positive random walk whose nonnegative i.i.d. increments $X_i$ satisfy the normalization condition $\mathbb{E}X_i = 1$. Note that if we proceed to the limit as $n \to \infty$, we consider a system of total mass one, which consists of, roughly speaking, infinitesimal particles homogeneously spread on $[0, 1]$; this is true for all the mentioned models of initial positions.

The mathematical interest in sticky particles systems arises mainly from relations between these systems and some nonlinear partial differential equations originating from fluid mechanics, for example, the Burgers equation. These equations admit interpretation in terms of sticky particles; see Gurbatov et al. [10], Brenier and Grenier [4] or E, Rykov and Sinai [6]. Sticky particles models are also used for numerical solving of other partial differential equations; see Chertock et al. [5] for explanations and further references.

As time goes, particles aggregate in clusters. Clusters become larger and larger while the number of clusters decreases until they merge into a single cluster containing all initial particles. This process of mass aggregation is strongly connected with additive coalescence; see Bertoin [2] and Giraud [9] for the most recent results and references.

The aggregation process resembles formation of a star from dispersed space dust and sticky particles models indeed have relations to astrophysics. It is appropriate to clarify these relations since they are not so direct and cause a lot of misunderstanding.

It is known that the distribution of galaxies in the universe is very inhomogeneous and the regions of high density form a peculiar cellular structure. The first attempt to understand the formation of such structures was made in 1970 by Zeldovich. Most of the mass in the universe is believed to exist in the form of particles that practically do not collide with each other and interact only gravitationally, for example, neutrinos. In his model, Zeldovich considered an initially homogeneous *collisionless* medium of particles moving by pure inertia; the gravitational interaction was taken away by an appropriate time change. He showed that singularities, that is, the thin regions of very high density of particles, so called "pancakes," appear even if initial speeds of particles form a smooth velocity field.

Zeldovich's approximate model, however, does not explain formation of the cellular structure of matter. His approximation does not take into account that particles hitting a "pancake" are hampered by its strong gravitational field and start oscillating inside the "pancake" instead of flying away. Although this gravitational adhesion of collisionless particles is not precisely the same as the real sticking, the model of sticky particles serves as a reasonable approximation. The effect of gravitational adhesion was then analyzed by the use of the Burgers equation; Gurbatov, Saichev and Shandarin proposed it in 1984 to extend Zeldovich's approximation, which is invalid after formation of "pancakes."

The model of sticky particles is *directly* mentioned in Gurbatov et al. [11]; a comprehensive survey of the formation of the Universe's large-scale structure could be found in Shandarin and Zeldovich [23].

1.2. *Statement of the problem and the results.* In general, the problem is to describe the process of mass aggregation. How fast is it? How large the clusters are? Where do clusters appear most intensively, and so forth? Numerous papers on the model (e.g., [8, 14, 16, 20, 25]) are dedicated to probabilistic description of various properties of the aggregation process as the number of initial particles $n$ tends to infinity. Thus, the behavior of a typical system consisting of a large number of particles is studied.

In this paper, we are interested in the asymptotic behavior of $K_n(t)$, which denotes the number of clusters at time $t$ in the system with $n$ initial particles. This variable is a decreasing random step function satisfying $K_n(0) = n$ and $K_n(t) = 1$ for $t \geq T_n^{\text{last}}$, where $T_n^{\text{last}}$ denotes the moment of the last collision. While calculating $K_n(t)$, we also count initial particles that have not experienced any collisions; in other words, $K_n(t)$ is the total number of particles existing at time $t$.

It is very important to know the behavior of $K_n(t)$. This gives us a deep understanding of the aggregation process since the average size of a cluster at time $t$ is $\frac{n}{K_n(t)}$.

At first we give a short deterministic example. Suppose that particles are located at points $\frac{1}{n}, \frac{2}{n}, \ldots, \frac{n}{n}$, that is, $S_i = i$. By simple calculations, we find that there would not be any collisions before $t = 1$. At the moment $t = 1$, all particles simultaneously stick together, hence $K_n(t) = n$ for $0 \leq t < 1$ and $K_n(t) = 1$ for $t \geq 1$.

However, when the initial positions are random, the aggregation process behaves entirely differently. In [25], the author proved the following statement.

FACT 1. *There exists a deterministic function $a(t)$ such that both in the Poisson and the uniform models of initial positions, for any $t \geq 0$, we have*

$$\frac{K_n(t)}{n} \xrightarrow{\mathbb{P}} a(t), \qquad n \to \infty. \tag{1}$$

*The function $a(t)$ is continuous, $a(0) = 1$, and $a(t) = 0$ for $t \geq 1$. We conjecture, on the basis of numerical simulations, that $a(t) = 1 - t^2$ for $0 \leq t \leq 1$.*

The relation $a(t) = 0$ for $t > 1$ is not of a surprise because we know from Giraud [8] that both in the Poisson and the uniform models, $T_n^{\text{last}} \xrightarrow{\mathbb{P}} 1$ (the limit constant is so "fine" due to the proper scaling of the model). Therefore, we say that the moment $t = 1$ is *critical*; note that this moment coincides with the moment of the total collision in the deterministic model.

The aim of this paper is to strengthen the result of [25]. We first generalize Fact 1 and prove it for the i.d. model. We will see [relations (19) and (27) below] that $a(t)$ is equal to the probability of a certain event that is expressed in terms of $X_i$. Also, we will prove that $a(t)$ depends on the common distribution of $X_i$ as follows: $a(t) = 1$ on $[0, \sqrt{\mu})$, where

$$\mu := \sup\{y : \mathbb{P}\{X_i < y\} = 0\};$$

$a(t) \in (0, 1)$ on $(\sqrt{\mu}, 1)$; and $a(t) = 0$ on $(1, \infty)$.

Furthermore, the recent results of the author [26] allow us to prove the conjecture from Fact 1 that $a^{\text{Poiss}}(t) = a^{\text{Unif}}(t) = 1 - t^2$ for $0 \leq t \leq 1$. There is an amazing contrast between the simplicity of this formula and the hard calculations one needs to obtain it. It is remarkable that now we know the limit function $a(t)$ for the main models of initial positions.

Our main goal is to improve (1) by finding the next term in the asymptotics of $K_n(t)$. The result is the following statement, where the standard symbol $\xrightarrow{\mathcal{D}}$ denotes weak convergence and $D$ denotes the Skorohod space.

THEOREM 1. *In the i.d. model with continuous $X_i$ satisfying $\mathbb{E}X_i^\gamma < \infty$ for some $\gamma > 4$, there exists a centered Gaussian process $K(\cdot)$ on $[0, 1)$ such that*

$$\frac{K_n(\cdot) - na(\cdot)}{\sqrt{n}} \xrightarrow{\mathcal{D}} K(\cdot) \qquad \text{in } D[0, 1 - \varepsilon] \text{ for all } \varepsilon \in (0, 1) \tag{2}$$

as $n \to \infty$. The process $K(\cdot)$ depends on the distribution of $X_i$. This process satisfies $K(0) = 0$ and has a.s. continuous trajectories. The covariance function $R(s,t)$ of $K(\cdot)$ is continuous on $[0,1)^2$, $R(s,t) > 0$ on $(\sqrt{\mu}, 1)^2$, and $R(s,t) = 0$ on $[0,1)^2 \setminus (\sqrt{\mu}, 1)^2$.

In the uniform model, (2) holds for some centered Gaussian process $K^{\text{Unif}}(\cdot)$ on $[0,1)$. This process satisfies $K^{\text{Unif}}(0) = 0$ and has a.s. continuous trajectories. The covariance function $R^{\text{Unif}}(s,t)$ of $K^{\text{Unif}}(\cdot)$ is continuous on $[0,1)^2$, and $R^{\text{Unif}}(s,t) = R^{\text{Poiss}}(s,t) - s^2 t^2$.

Thus, the Poisson and the uniform models lead to different limit processes $K^{\text{Poiss}}(\cdot)$ and $K^{\text{Unif}}(\cdot)$, although $a^{\text{Poiss}}(\cdot) = a^{\text{Unif}}(\cdot)$.

As an immediate corollary of Theorem 1 (see Billingsley [3], Section 15), we get

$$(3) \qquad \frac{K_n(t) - na(t)}{\sqrt{n}} \xrightarrow{\mathcal{D}} \mathcal{N}(0, \sigma^2(t)), \qquad n \to \infty$$

for any $t < 1$, where $\sigma^2(t) := R(t,t)$. It is possible to show that in the i.d. model, (3) holds for all $t \neq 1$ under the less restrictive condition $\mathbb{E} X_i^2 < \infty$, with $\sigma^2(t) = 0$ for $t > 1$; continuity of $X_i$ is not required.

We also study convergence of the left-hand side of (3) at the critical moment $t = 1$. Apparently, the limit is not Gaussian, but this complicated problem is related to a curious, but hardly provable conjecture on integrated random walks. In view of this non-Gaussianity, it seems impossible to prove any extended version of Theorem 1 that describes the weak convergence of trajectories on the whole interval $[0,1]$; we refer to Section 7 for further discussion.

We finish this subsection with a note on scaling. In our model, the masses of particles are equal to $\frac{1}{n}$ and the distances between them are of the order $\frac{1}{n}$. Let us rescale the i.d. model by multiplying all masses and distances by $n$: the system of particles of mass one each, initially located at points $S_1 - S_{[n/2]}, S_2 - S_{[n/2]}, \ldots, S_n - S_{[n/2]}$, is called the *expanding* model. The particles are shifted by $S_{[n/2]}$ because we want the system to expand "filling" the whole line as $n \to \infty$ rather than only the positive half-line.

All results of our paper hold true for the expanding model. This is not unexpected because the shift does not produce any changes and the rescaling of masses is equivalent to the time contraction by $n$ times while the rescaling of distances is equivalent to the time expansion by $n$ times. We refer the reader to Section 2 below or to Lifshits and Shi [16] for rigorous arguments.

1.3. *Organization of the paper.* In Section 2 we describe a general method which is used to study systems of sticky particles. This method is applied for studying the i.d. model in Section 3, where we investigate some properties of the aggregation process. We will show that *the aggregation process is highly local*, that is, the behavior of a particle is essentially defined by the motion of neighbor particles. This localization property suggests that we could use limit theorems for weakly dependent variables to prove both Fact 1 and Theorem 1 for the i.d. model; this will be done in Section 4. Then we will prove Theorem 1 for the uniform model in Section 5. In Section 6 we study the number of clusters at the critical moment $t = 1$. Some open questions are discussed in Section 7.

**2. Method of barycenters.** In this section we briefly describe the method of barycenters, which is the main tool used to study systems of sticky particles; it is also applicable to more general models where particles could have nonzero initial speeds and different masses. The method of barycenters was independently introduced by E, Rykov and Sinai [6] and Martin and Piasecki [20].

Let us start with several definitions. We always numerate particles from left to right and identify particles with their numbers. A *block* of particles is a nonempty set $J \subset [1, n]$ consisting of consecutive numbers. For example, the block $(i, i+k]$ consists of particles $i+1, \ldots, i+k$. Note that there are not any relations between blocks and clusters: for example, a block's particles could be contained in different clusters and these clusters could even contain particles that do not belong to the block.

It is convenient to assume that initial particles do not vanish at collisions but continue to exist in created clusters. Then the coordinate $x_{i,n}(t)$ of a particle $i$ could be defined as the coordinate of a cluster that contains the particle at time $t$. The second subscript $n$ always indicates the number of initial particles; we will omit this subscript as often as possible.

By $x_J(t) := |J|^{-1} \sum_{i \in J} x_i(t)$ denote the position of the *barycenter* of a block $J$ at time $t$. Further, define

$$x_J^*(t) := x_J(0) + \tfrac{1}{2}(M_J^{(R)} - M_J^{(L)})t^2,$$

where $M_J^{(R)} := n^{-1}(n - \max_{j \in J})$ and $M_J^{(L)} := n^{-1}(\min_{j \in J} - 1)$ are the total masses of particles located to the right and to the left of the block $J$, respectively.

A block is *free from the right* up to time $t$ if, up to this time, the block's particles did not collide with particles initially located to the right of the block. We similarly define blocks that are *free from the left* and say that a block is *free* up to time $t$ if it is both free from the right and from the left.

The next statement plays the key role in the analysis of sticky particles systems. *The barycenter of a free block moves as an imaginary particle consisting of all particles of the block put together at the initial barycenter.* In a more precise and general way, we state the following.

PROPOSITION 1. *If a block $J$ is free from the right (resp. left) up to time $t$, then $x_J(s) \geq x_J^*(s)$ for $s \in [0,t]$ [resp. $x_J(s) \leq x_J^*(s)$]. If a block $J$ is free up to time $t$, then $x_J(s) = x_J^*(s)$ for $s \in [0,t]$.*

This statement could be found, for example, in Lifshits and Shi [16], Proposition 4.1. The easy proof is based on the property of conservation of momentum.

The moment when a particle $j$ sticks with its right-hand side neighbor $j+1$ is called the *merging time* $T_{j,n}$ of the particle $j$. In other words, $T_{j,n}$ is the first moment when particles $j$ and $j+1$ are contained in a common cluster; here $j \in [1, n-1]$. Proposition 4.3 from Lifshits and Shi [16], which is stated below, gives us a way to calculate $T_{j,n}$.

PROPOSITION 2. *For every $j \in [1, n-1]$, we have*

$$(4) \qquad T_{j,n} = \min_{\substack{j < k \leq n \\ 0 \leq l < j}} \{s \geq 0 : x_{(j,k]}^*(s) = x_{(l,j]}^*(s)\}.$$

Thus, $T_{j,n}$ is expressed by means of barycenters. Note that since

$$(5) \qquad x_{(j,k]}^*(s) - x_{(l,j]}^*(s) = x_{(j,k]}(0) - x_{(l,j]}(0) - \frac{k-l}{2n}s^2,$$

each of the equations $x_{(j,k]}^*(s) = x_{(l,j]}^*(s)$ has a unique nonnegative solution. We also mention that at the moment $T_{j,n}$ appears a cluster that consists of the particles $l+1, \ldots, k$, where $k$ and $l$ are minimizers of the right-hand side of (4).

We will prove Proposition 2 since the proof is simple and perfectly illustrates the sense of the method of barycenters.

PROOF OF PROPOSITION 2. For any $u < T_{j,n}$, the particles $j$ and $j+1$ are contained in different clusters. Therefore, for every $l < j$, the block $[l,j]$ is free from the right up to time $u$, and for every $k > j$, the block $[j+1,k]$ is free from the left. By Proposition 1,

$$x_{(l,j]}^*(u) \leq x_{(l,j]}(u) \leq x_j(u) < x_{j+1}(u)$$
$$\leq x_{(j,k]}(u) \leq x_{(j,k]}^*(u),$$

and since, by (5), the function $x_{(j,k]}^*(s) - x_{(l,j]}^*(s)$ is decreasing for $s \geq 0$, we conclude that

$$u < \{s \geq 0 : x_{(j,k]}^*(s) = x_{(l,j]}^*(s)\}.$$

Taking minimum over $k, l$ and taking supremum over $u$, we get $T_{j,n} \leq \min\{\cdots\}$.

Let us prove the last inequality in the other direction. By the definition of $T_{j,n}$, there exist an $l < j$ and a $k > j$ such that the blocks $(l,j]$ and $(j,k]$ are free up to time $T_{j,n}$ (clusters containing particles from these blocks collide exactly at time $T_{j,n}$). In view of Proposition 1,

$$x^*_{(l,j]}(T_{j,n}) = x_{(l,j]}(T_{j,n}) = x_{(j,k]}(T_{j,n}) = x^*_{(j,k]}(T_{j,n});$$

hence $T_{j,n} = \{s \geq 0 : x^*_{(j,k]}(s) = x^*_{(l,j]}(s)\}$ and $T_{j,n} \geq \min\{\cdots\}$. □

**3. Study of the i.d. model. The localization property.** At first, note that

$$(6) \qquad K_n(t) = 1 + \sum_{i=1}^{n-1} \mathbb{1}_{\{t < T_{i,n}\}}$$

because the total number of clusters decreases by one at each moment $T_{i,n}$. This representation plays the key role in the investigation of $K_n(t)$. Clearly, we need to study properties of the r.v.'s $T_{i,n}$ to prove limit theorems for $K_n(t)$; such study will be done in this section.

3.1. *The initial study.* Let us simplify the representation for $T_{j,n}$ from Proposition 2. In this section we consider the i.d. model of initial positions, where $x_{j,n}(0) = \frac{1}{n} S_j$. Recall that $S_j$ is a random walk with i.i.d. increments $\{X_j\}_{j \in \mathbb{Z}}$ (we will need the variables $\{X_j\}_{j \leq 0}$ later).

Rewrite the initial distance between barycenters as

$$x_{(j,k]}(0) - x_{(l,j]}(0)$$

$$= \frac{1}{k-j} \sum_{i=j+1}^{k} \frac{1}{n} S_i - \frac{1}{j-l} \sum_{i=l+1}^{j} \frac{1}{n} S_i$$

$$= \frac{1}{n}\left(\frac{1}{k-j} \sum_{i=j+1}^{k} (S_i - S_{j+1}) + \frac{1}{j-l} \sum_{i=l+1}^{j}(S_j - S_i) + (S_{j+1} - S_j)\right)$$

$$= \frac{1}{n}\left(\frac{1}{k-j} \sum_{i=1}^{k-j-1}(S_{j+i+1} - S_{j+1}) + \frac{1}{j-l}\sum_{i=1}^{j-l-1}(S_j - S_{j-i}) + X_{j+1}\right);$$

let us agree that $\sum_{\varnothing} := 0$. Further, by

$$x_{(j,k]}(0) - x_{(l,j]}(0)$$

$$= \frac{1}{n}\left(\frac{1}{k-j}\sum_{i=1}^{k-j-1}\sum_{m=j+2}^{j+i+1} X_m + \frac{1}{j-l}\sum_{i=1}^{j-l-1}\sum_{m=j-i+1}^{j} X_m + X_{j+1}\right)$$

$$= \frac{1}{n}\left(\frac{1}{k-j}\sum_{i=1}^{k-j-1}(k-j-i)X_{j+i+1}\right.$$

$$+ \frac{1}{j-l}\sum_{i=1}^{j-l-1}(j-l-i)X_{j-i+1} + X_{j+1}\bigg),$$

and (5), we have

$$x^*_{(j,k]}(s) - x^*_{(l,j]}(s) = \frac{1}{n}F_{k-j,j,j-l}(s),$$

where

(7)
$$F_{p,j,q}(s) := \frac{1}{p}\sum_{i=1}^{p-1}(p-i)X_{j+i+1}$$
$$+ \frac{1}{q}\sum_{i=1}^{q-1}(q-i)X_{j-i+1} + X_{j+1} - \frac{p+q}{2}s^2$$

(for $p, q \geq 1$ and $j \in \mathbb{Z}$). Now, by Proposition 2, we get

$$T_{j,n} = \min_{\substack{j<k\leq n \\ 0\leq l<j}}\{s\geq 0 : F_{k-j,j,j-l}(s) = 0\}$$

(8)
$$= \min_{\substack{1\leq k\leq n-j \\ 1\leq l\leq j}}\{s\geq 0 : F_{k,j,l}(s) = 0\}.$$

Note that $F_{p,j,q}(0) \geq 0$ for all $p, j, q$ and $F_{p,j,q}(s)$ is decreasing for $s \geq 0$. This function could be also written in the more convenient form:

(9)
$$F_{p,j,q}(s) = \frac{1}{p}\sum_{i=1}^{p-1}(p-i)(X_{j+i+1} - s^2)$$
$$+ \frac{1}{q}\sum_{i=1}^{q-1}(q-i)(X_{j-i+1} - s^2) + (X_{j+1} - s^2).$$

3.2. *Localization property of the aggregation process.* We see that $T_{j,n}$ is a function of $X_2, \ldots, X_n$; in other words, it is necessary to know the distances between all $n$ particles to find $T_{j,n}$. The aggregation process is actually highly local, that is, *the value of $T_{j,n}$ is essentially defined by the initial distances between neighbor particles $\{i\}$ of $j$ for which $|j - i|$ is small enough.*

To make this statement rigorous, we need to introduce the following notation. Let us put

$$T_j^{(M)} := \min_{1\leq k,l\leq M}\{s \geq 0 : F_{k,j,l}(s) = 0\}, \qquad j \in \mathbb{Z}, M \in \mathbb{N},$$

which is expressed in terms of the variables $\{X_i\}_{|j-i|\leq M}$ only. Also, define

$$T_j := \inf_{k,l\geq 1}\{s \geq 0 : F_{k,j,l}(s) = 0\}, \qquad j \in \mathbb{Z},$$

which is, in some sense, the merging time in an appropriate infinite system of particles. The reader could construct such system by considering the limit of the expanding model, see Section 1.

It is clear that

$$T_j \leq T_{j,n} \leq T_j^{(j \wedge n-j)}, \qquad j, n \in \mathbb{N}, j \leq n, \tag{10}$$

where by $\wedge$ and $\vee$ we denote minimum and maximum, respectively, and

$$T_j \leq T_j^{(M)}, \qquad j \in \mathbb{Z}, M \in \mathbb{N}. \tag{11}$$

Let us estimate the rate of the convergence of $\mathbb{P}\{T_j \neq T_j^{(M)}\}$ to zero as the "radius of the neighborhood" $M$ tends to infinity. We thus could "measure" the above-mentioned locality of the aggregation process. In fact, by (10), we have $\mathbb{P}\{T_{j,n} \neq T_j^{(M)}\} \leq \mathbb{P}\{T_j \neq T_j^{(M)}\}$ for any $n \in \mathbb{N}$, $j \leq n$, and $M \leq j \wedge n-j$.

LEMMA 1. *Suppose $\mathbb{E}X_i^\gamma < \infty$ for some $\gamma \geq 1$. Then there exists a non-decreasing function $\rho(t)$ such that*

$$\max(\mathbb{P}\{\mathbb{1}_{\{t \leq T_j\}} \neq \mathbb{1}_{\{t \leq T_j^{(M)}\}}\}, \mathbb{P}\{T_j \neq T_j^{(M)}, T_j^{(M)} \leq t\}) \leq \rho(t) M^{1-\gamma} \tag{12}$$

*for any $t \in (0,1)$, $j \in \mathbb{Z}$, and $M \in \mathbb{N}$. Moreover, for any $t < 1$, the left-hand side of (12) is $o(M^{1-\gamma})$.*

PROOF. Let us estimate the first probability in the left-hand side of (12). By properties of $F_{k,j,l}(\cdot)$ and definitions of $T_j^{(M)}$ and of $T_j$,

$$\mathbb{P}\{\mathbb{1}_{\{t \leq T_j\}} \neq \mathbb{1}_{\{t \leq T_j^{(M)}\}}\} = \mathbb{P}\{T_j < t \leq T_j^{(M)}\}$$

$$= \mathbb{P}\Big\{\inf_{k,l \geq 1} F_{k,j,l}(t) < 0, \min_{1 \leq k,l \leq M} F_{k,j,l}(t) \geq 0\Big\}.$$

By (9), this expression does not depend on $j$, and putting $j := -1$,

$$\mathbb{P}\{\mathbb{1}_{\{t \leq T_j\}} \neq \mathbb{1}_{\{t \leq T_j^{(M)}\}}\}$$

$$= \mathbb{P}\Big\{\inf_{k \geq 1} \frac{1}{k}\sum_{i=1}^{k-1}(k-i)(X_i - t^2)$$

$$+ \inf_{l \geq 1} \frac{1}{l}\sum_{i=1}^{l-1}(l-i)(X_{-i} - t^2) + (X_0 - t^2) < 0,$$

$$\min_{1 \leq k \leq M} \frac{1}{k}\sum_{i=1}^{k-1}(k-i)(X_i - t^2)$$

CLUSTERING IN A STOCHASTIC MODEL OF ONE-DIMENSIONAL GAS 11$$+ \min_{1 \leq l \leq M} \frac{1}{l} \sum_{i=1}^{l-1}(l-i)(X_{-i} - t^2) + (X_0 - t^2) \geq 0 \Big\}.$$

We then compare the inequalities in the braces and obtain

$$\mathbb{P}\{\mathbb{1}_{\{t \leq T_j\}} \neq \mathbb{1}_{\{t \leq T_j^{(M)}\}}\}$$

$$\leq 2\mathbb{P}\Big\{\inf_{k>M} \frac{1}{k}\sum_{i=1}^{k-1}(k-i)(X_i - t^2) < \min_{1 \leq k \leq M} \frac{1}{k}\sum_{i=1}^{k-1}(k-i)(X_i - t^2)\Big\}$$

$$= 2\mathbb{P}\Big\{\inf_{k>M} \frac{1}{k}\sum_{i=1}^{k-1}(S_i - it^2) < \min_{1 \leq k \leq M} \frac{1}{k}\sum_{i=1}^{k-1}(S_i - it^2)\Big\}$$

$$\leq 2\mathbb{P}\Big\{\inf_{k>M} \frac{1}{k}\sum_{i=1}^{k-1}(S_i - it^2) < \min_{k \in \{1,M\}} \frac{1}{k}\sum_{i=1}^{k-1}(S_i - it^2)\Big\}.$$

Now rewrite the event in the last line as

$$\Big\{\exists k > M : \frac{1}{k}\sum_{i=1}^{k-1}(S_i - it^2) < \min\Big(0, \frac{1}{M}\sum_{i=1}^{M-1}(S_i - it^2)\Big)\Big\}$$

$$= \Big\{\exists k > M : \frac{1}{k}\sum_{i=1}^{M-1}(S_i - it^2)$$

$$+ \frac{1}{k}\sum_{i=M}^{k-1}(S_i - it^2) < \min\Big(0, \frac{1}{M}\sum_{i=1}^{M-1}(S_i - it^2)\Big)\Big\}.$$

Analyzing both cases $0 \leq \frac{1}{M}\sum_{i=1}^{M-1}(S_i - it^2)$ and $0 > \frac{1}{M}\sum_{i=1}^{M-1}(S_i - it^2)$, we conclude that the considered event implies

$$\Big\{\exists k > M : \frac{1}{k}\sum_{i=M}^{k-1}(S_i - it^2) < 0\Big\} = \Big\{\exists k > M : \sum_{i=M}^{k-1}(S_i - it^2) < 0\Big\}.$$

Clearly, the latter implies

$$\{\exists i \geq M : S_i - it^2 < 0\} = \Big\{\inf_{i \geq M} \frac{S_i}{i} < t^2\Big\};$$

hence, combining all the estimates together, we get

(13) $$\mathbb{P}\{\mathbb{1}_{\{t \leq T_j\}} \neq \mathbb{1}_{\{t \leq T_j^{(M)}\}}\} \leq 2\mathbb{P}\Big\{\inf_{i \geq M} \frac{S_i}{i} < t^2\Big\}.$$

Note that we obtained (13) without any assumptions on the moments of $X_i$.

We now estimate the right-hand side of (13); recall that $\mathbb{E}X_i = 1$. Then the first part of (12) immediately follows from the classical result of Baum and Katz [1] (see their Theorem 3 and Lemma):

FACT 2. *If $\mathbb{E}X_i = a$ and $\mathbb{E}|X_i|^\gamma < \infty$ for some $\gamma \geq 1$, then*

$$\mathbb{P}\left\{\sup_{i \geq k}\left|\frac{S_i}{i} - a\right| > \varepsilon\right\} = o(k^{1-\gamma}), \qquad k \to \infty$$

*for any $\varepsilon > 0$. In addition, the series $\sum_{k=1}^\infty \mathbb{P}\{\sup_{i \geq k} |\frac{S_i}{i} - a| > \varepsilon\}$ converges for all $\varepsilon > 0$ if $\gamma = 2$.*

The estimation of the second probability in the left-hand side of (12) is completely analogous, since

$$\{T_j \neq T_j^{(M)}, T_j^{(M)} \leq t\}$$
$$= \{T_j < T_j^{(M)} \leq t\}$$
$$= \left\{\inf_{1 \leq k, l} F_{k,j,l}(T_j^{(M)}) < 0, \min_{1 \leq k, l \leq M} F_{k,j,l}(T_j^{(M)}) = 0, T_j^{(M)} \leq t\right\}.$$

We put $j := -1$, repeat the estimates, and get

$$\mathbb{P}\{T_j \neq T_j^{(M)}, T_j^{(M)} \leq t\} \leq 2\mathbb{P}\{\exists i \geq M : S_i - i[T_{-1}^{(M)}]^2 < 0, T_{-1}^{(M)} \leq t\}$$

instead of (13). The right-hand side does not exceed $2\mathbb{P}\{\exists i \geq M : S_i - it^2 < 0\}$, hence

$$(14) \qquad \mathbb{P}\{T_j \neq T_j^{(M)}, T_j^{(M)} \leq t\} \leq 2\mathbb{P}\left\{\inf_{i \geq M} \frac{S_i}{i} < t^2\right\}. \qquad \square$$

3.3. *The distribution function of $T_0$ in the Poisson model.* It is amazing that in the Poisson model, the distribution function of $T_0$ could be found explicitly. This is important because by (27) below, the limit function $a(t)$ equals $\mathbb{P}\{T_0 > t\}$ for the i.d. model. Also, in the proof of Theorem 1 for the uniform model, we will need $a^{\text{Poiss}}(t) = \mathbb{P}\{T_0^{\text{Poiss}} \geq t\}$ to be twice differentiable and have a continuous second derivative.

LEMMA 2. *In the Poisson model, for $0 \leq t \leq 1$, we have*

$$(15) \qquad \mathbb{P}\{T_0 \geq t\} = 1 - t^2.$$

*In addition, for $t \geq 0$, $n \geq 2$, and $1 \leq j \leq n-1$, we have*

$$(16) \quad \mathbb{P}\{T_{j,n} \geq t\} = e^{t^2} \mathbb{P}\left\{\min_{1 \leq k \leq j} \sum_{i=1}^k (S_i - it^2) \geq 0\right\}$$
$$\times \mathbb{P}\left\{\min_{1 \leq k \leq n-j} \sum_{i=1}^k (S_i - it^2) \geq 0\right\},$$

*where $S_i$ is a standard exponential random walk.*

PROOF. We start with (16). By (8), (9) and properties of $F_{k,j,l}(\cdot)$,

$$\mathbb{P}\{T_{j,n} \geq t\} = \mathbb{P}\left\{\min_{\substack{1 \leq k \leq n-j \\ 1 \leq l \leq j}} F_{k,j,l}(t) \geq 0\right\}$$

(17)
$$= \mathbb{P}\left\{\min_{1 \leq k \leq n-j} \frac{1}{k} \sum_{i=1}^{k-1}(k-i)(X_{j+i+1} - t^2) \right.$$
$$\left. + \min_{1 \leq l \leq j} \frac{1}{l} \sum_{i=1}^{l-1}(l-i)(X_{j-i+1} - t^2) + X_{j+1} - t^2 \geq 0\right\}.$$

In the right-hand side of the last equality, by $Y$ denote the first minimum and by $\tilde{Y}$ denote the second one.

Suppose $X$ is a standard exponential r.v., $Z$ is a nonnegative r.v., and that $X$ and $Z$ are independent; then

$$\mathbb{P}\{Z \leq X\} = \int_0^\infty \mathbb{P}\{Z \leq x\} e^{-x}\, dx$$
$$= \int_0^\infty \mathbb{E}\mathbb{1}_{\{Z \leq x\}} e^{-x}\, dx = \mathbb{E}\int_0^\infty \mathbb{1}_{\{Z \leq x\}} e^{-x}\, dx = \mathbb{E}e^{-Z}.$$

Hence in view of independence of $Y$, $\tilde{Y}$, $X_{j+1}$ we get

$$\mathbb{P}\{Y + \tilde{Y} + X_{j+1} - t^2 \geq 0\} = \mathbb{E}e^{Y + \tilde{Y} - t^2} = e^{t^2}\mathbb{E}e^{Y-t^2}\mathbb{E}e^{\tilde{Y}-t^2};$$

and therefore,

$$\mathbb{P}\{T_{j,n} \geq t\} = e^{t^2}\mathbb{P}\{Y + X_{j+1} - t^2 \geq 0\} \cdot \mathbb{P}\{\tilde{Y} + X_{j+1} - t^2 \geq 0\}.$$

Now, by

$$\mathbb{P}\{\tilde{Y} + X_{j+1} - t^2 \geq 0\}$$
(18)
$$= \mathbb{P}\left\{\min_{1 \leq l \leq j} \frac{1}{l} \sum_{i=1}^{l-1}(l-i)(X_{j-i+1} - t^2) + X_{j+1} - t^2 \geq 0\right\}$$
$$= \mathbb{P}\left\{\min_{1 \leq l \leq j} \left(\sum_{i=1}^{l-1}(l-i)(X_{i+1} - t^2) + l(X_1 - t^2)\right) \geq 0\right\}$$
$$= \mathbb{P}\left\{\min_{1 \leq l \leq j} \sum_{i=1}^{l}(l-i+1)(X_i - t^2) \geq 0\right\},$$

we conclude the proof of (16). Indeed, the expression in the last line equals the first probability in the right-hand side of (16).

Now let us prove (15). From the definition of $T_0$ and $T_0^{(k)}$ we see that $\mathbb{1}_{\{t \leq T_0^{(k)}\}} \to \mathbb{1}_{\{t \leq T_0\}}$ a.s. as $k \to \infty$; then by (16),

$$\mathbb{P}\{T_0 \geq t\} = e^{t^2} \mathbb{P}^2 \left\{ \inf_{k \geq 1} \sum_{i=1}^{k} (S_i - it^2) \geq 0 \right\}.$$

Then we need to check that

$$\mathbb{P}\left\{ \inf_{k \geq 1} \sum_{i=1}^{k} (S_i - it) \geq 0 \right\} = \sqrt{1-t}\, e^{-t/2}$$

for $0 \leq t \leq 1$. The complicated calculations of this probability take more then ten pages. Therefore, they were separated into independent paper [26]. Although these calculations seem to be technical, they are based on quite original ideas. □

3.4. *Some properties of the variables $T_i$.* In this subsection we prove several important properties of the r.v.'s $T_i$.

1. The sequence $T_i$ is stationary.

*Proof.* This statement immediately follows from the definition of $T_i$ and stationarity of $X_i$, which are i.i.d.

2. The common distribution function of $T_i$ is defined by

$$\begin{aligned}
\mathbb{P}\{T_i \geq t\} = \mathbb{P}\Bigg\{ & \inf_{k \geq 1} \frac{1}{k} \sum_{i=1}^{k-1} (k-i)(X_i - t^2) \\
& + \inf_{l \geq 1} \frac{1}{l} \sum_{i=1}^{l-1} (l-i)(X_{-i} - t^2) + (X_0 - t^2) \geq 0 \Bigg\}.
\end{aligned} \tag{19}$$

*Proof.* This formula follows from (9).

3. We have $\mathbb{P}\{\sqrt{\mu} \leq T_i \leq 1\} = 1$ while $\sup\{y : \mathbb{P}\{T_i < y\} = 0\} = \sqrt{\mu}$ and $\inf\{y : \mathbb{P}\{T_i < y\} = 1\} = 1$; recall that $\mu = \sup\{y : \mathbb{P}\{X_i < y\} = 0\}$. In addition, if $0 < \mathbb{D}X_i < \infty$, then $\mathbb{P}\{T_i = 1\} = 0$.

*Proof.* First, $\mathbb{P}\{\sqrt{\mu} \leq T_i\} = 1$ is trivial, because both infima in (19) are nonpositive.

Second, fix a $t \geq 1$ and consider $\mathbb{P}\{T_i \geq t\}$. Taking into account that infima in (19) are nonpositive, we obtain

$$\mathbb{P}\{T_i \geq t\} \leq \mathbb{P}\left\{ \inf_{k \geq 1} \frac{1}{k} \sum_{i=1}^{k-1} (k-i)(X_i - t^2) + (X_0 - t^2) \geq 0 \right\}.$$

Then by the same arguments as in (18),

$$\mathbb{P}\{T_i \geq t\} \leq \mathbb{P}\left\{ \inf_{k \geq 1} \sum_{i=1}^{k} (k-i+1)(X_i - t^2) \geq 0 \right\} = \mathbb{P}\left\{ \inf_{k \geq 1} \sum_{i=1}^{k} (S_i - it^2) \geq 0 \right\}.$$

By the strong law of large numbers, this probability is zero for all $t > 1$.

If $t = 1$ and $0 < \mathbb{D}X_i < \infty$, then

$$\mathbb{P}\left\{\inf_{k \geq 1} \sum_{i=1}^{k}(S_i - i) \geq 0\right\} = \lim_{n \to \infty} \mathbb{P}\left\{\min_{1 \leq k \leq n} \sum_{i=1}^{k}(S_i - i) \geq 0\right\}$$

$$= \lim_{n \to \infty} \mathbb{P}\left\{\min_{1 \leq k \leq n} \frac{1}{n} \sum_{i=1}^{k} \frac{S_i - i}{\sqrt{n\mathbb{D}X_i}} \geq 0\right\},$$

and from the invariance principle, we get

$$\mathbb{P}\{T_i \geq 1\} \leq \mathbb{P}\left\{\min_{0 \leq s \leq 1} \int_0^s W(u)\,du \geq 0\right\}.$$

It follows from the asymptotics of unilateral small deviation probabilities of an integrated Wiener process, see (43) and (44) below, that the last expression equals zero.

Third, $\sup\{y : \mathbb{P}\{T_i < y\} = 0\} = \sqrt{\mu}$ and $\inf\{y : \mathbb{P}\{T_i < y\} = 1\} = 1$ follow if we prove that for any $t < \mathbb{E}X_i = 1$, the common distribution of the i.i.d. infima in (19) has an atom at zero. But we have

$$\mathbb{P}\left\{\inf_{k \geq 1} \frac{1}{k} \sum_{i=1}^{k-1}(k-i)(X_i - t^2) = 0\right\} = \mathbb{P}\left\{\inf_{k \geq 1} \frac{1}{k} \sum_{i=1}^{k-1}(S_i - it^2) = 0\right\}$$

$$\geq \mathbb{P}\left\{\inf_{i \geq 1} \frac{S_i}{i} \geq t^2\right\},$$

and it could be shown via the strong law of large numbers that the last probability is strictly positive for all $t < 1$.

4. Suppose $X_i$ is continuous. Then $T_j^{(k)}$ and $T_{j,n}$ are continuous for any $j, k, n$ and the common distribution of $T_j$ could have an atom only at 1. In addition, if $\mathbb{E}X_i^2 < \infty$, then $T_j$ are continuous.

*Proof.* By (7) and (8),

(20) $$T_{j,n} = \min_{\substack{1 \leq k \leq n-j \\ 1 \leq l \leq j}} \sqrt{H(k, j, l)},$$

where

$$H(p, j, q) := \frac{2}{p+q}\left(\frac{1}{p}\sum_{i=1}^{p-1}(p-i)X_{j+i+1} + \frac{1}{q}\sum_{i=1}^{q-1}(q-i)X_{j-i+1} + X_{j+1}\right).$$

Hence $T_{j,n}$ is continuous as a minimum of a finite number of continuous r.v.'s. The $T_j^{(k)}$ are also continuous because $T_j^{(k)} \stackrel{\mathcal{D}}{=} T_{k,2k}$.

Now we prove the continuity of $T_j$. By Property 3, it only remains to verify that $\mathbb{P}\{T_j \geq t\}$ is continuous on $[0,1)$. But $\mathbb{P}\{T_j^{(k)} \geq t\} - \mathbb{P}\{T_j \geq t\} = \mathbb{P}\{\mathbb{1}_{\{t \leq T_j\}} \neq \mathbb{1}_{\{t \leq T_j^{(k)}\}}\}$, and in view of (13),

$$\sup_{0 \leq t \leq s} |\mathbb{P}\{T_j^{(k)} \geq t\} - \mathbb{P}\{T_j \geq t\}| \leq \sup_{0 \leq t \leq s} 2\mathbb{P}\left\{\inf_{m \geq k} \frac{S_m}{m} < t^2\right\} = 2\mathbb{P}\left\{\inf_{m \geq k} \frac{S_m}{m} < s^2\right\}$$

for every $s < 1 = \mathbb{E}X_i$. The last expression tends to zero by the strong law of large numbers; then $\mathbb{P}\{T_j \geq t\}$ is continuous on $[0,s]$ as a uniform limit of continuous functions $\mathbb{P}\{T_j^{(k)} \geq t\}$. Since $s < 1$ is arbitrary, $\mathbb{P}\{T_j \geq t\}$ is continuous on $[0,1)$.

5. The $\text{cov}(\mathbb{1}_{\{s \leq T_0\}}, \mathbb{1}_{\{t \leq T_k\}})$ tends to zero as $k \to \infty$ for all $s, t \in [0,1)$. If, in addition, $\mathbb{E}X_i^\gamma < \infty$ for some $\gamma > 1$, then for any $s, t \in [0,1)$ and $k \in \mathbb{N}$, we have

(21) $$|\text{cov}(\mathbb{1}_{\{s \leq T_0\}}, \mathbb{1}_{\{t \leq T_k\}})| \leq 2^\gamma (\rho(s) + \rho(t))k^{1-\gamma}.$$

*Proof.* The idea is to approximate $\mathbb{1}_{\{s \leq T_0\}}$ and $\mathbb{1}_{\{t \leq T_k\}}$ by $\mathbb{1}_{\{s \leq T_0^{(k/2)}\}}$ and $\mathbb{1}_{\{t \leq T_k^{(k/2)}\}}$, respectively; here by $k/2$ we mean $\lceil k/2 \rceil$, where $\lceil x \rceil = \min\{m \in \mathbb{Z} : m \geq x\}$. Note that $\mathbb{1}_{\{s \leq T_0^{(k/2)}\}}$ and $\mathbb{1}_{\{t \leq T_k^{(k/2)}\}}$ are independent because the first is a function of $\{X_i\}_{i \leq k/2}$ while the second is a function of $\{X_i\}_{i \geq k/2+1}$. We then have

$$|\text{cov}(\mathbb{1}_{\{s \leq T_0\}}, \mathbb{1}_{\{t \leq T_k\}})|$$
$$= |\text{cov}(\mathbb{1}_{\{s \leq T_0\}}, \mathbb{1}_{\{t \leq T_k\}}) - \text{cov}(\mathbb{1}_{\{s \leq T_0^{(k/2)}\}}, \mathbb{1}_{\{t \leq T_k^{(k/2)}\}})|$$
$$\leq |\mathbb{E}(\mathbb{1}_{\{s \leq T_0\}}\mathbb{1}_{\{t \leq T_k\}} - \mathbb{1}_{\{s \leq T_0^{(k/2)}\}}\mathbb{1}_{\{t \leq T_k^{(k/2)}\}})|$$
(22) $$\quad + |\mathbb{E}(\mathbb{1}_{\{s \leq T_0\}} - \mathbb{1}_{\{s \leq T_0^{(k/2)}\}})| + |\mathbb{E}(\mathbb{1}_{\{t \leq T_k\}} - \mathbb{1}_{\{t \leq T_k^{(k/2)}\}})|$$
$$= \mathbb{P}\{\mathbb{1}_{\{s \leq T_0\}}\mathbb{1}_{\{t \leq T_k\}} \neq \mathbb{1}_{\{s \leq T_0^{(k/2)}\}}\mathbb{1}_{\{t \leq T_k^{(k/2)}\}}\}$$
$$\quad + \mathbb{P}\{\mathbb{1}_{\{s \leq T_0\}} \neq \mathbb{1}_{\{s \leq T_0^{(k/2)}\}}\} + \mathbb{P}\{\mathbb{1}_{\{t \leq T_k\}} \neq \mathbb{1}_{\{t \leq T_k^{(k/2)}\}}\}.$$

But

$$\mathbb{P}\{\mathbb{1}_{\{s \leq T_0\}}\mathbb{1}_{\{t \leq T_k\}} \neq \mathbb{1}_{\{s \leq T_0^{(k/2)}\}}\mathbb{1}_{\{t \leq T_k^{(k/2)}\}}\}$$
$$\leq \mathbb{P}\{\mathbb{1}_{\{s \leq T_0\}} \neq \mathbb{1}_{\{s \leq T_0^{(k/2)}\}} \cup \mathbb{1}_{\{t \leq T_k\}} \neq \mathbb{1}_{\{t \leq T_k^{(k/2)}\}}\},$$

therefore the result follows from Lemma 1.

6. The r.v.'s $\{T_i\}_{i \in \mathbb{Z}}$, $\{T_i^{(k)}\}_{i \in \mathbb{Z}}$, and $\{T_{i,n}\}_{i=1}^{n-1}$ are *associated*; the author owes this observation to M. A. Lifshits.

*Proof.* Let us first recall the definition and some basic properties of associated variables. R.v.'s $\xi_1, \ldots, \xi_m$ are *associated* if for any coordinate-wise nondecreasing functions $f, g : \mathbb{R}^m \to \mathbb{R}$, it is true that

$$\mathrm{cov}(f(\xi_1, \ldots, \xi_m), g(\xi_1, \ldots, \xi_m)) \geq 0$$

(assuming that the left-hand side is well defined). An infinite set of r.v.'s is associated if any finite subset of its variables is associated.

The following sufficient conditions of association are well known; see [7].

(a) Independent variables are associated.

(b) Coordinate-wise nondecreasing functions (of finite number of arguments) of associated r.v.'s are associated.

(c) If the variables $\xi_{1,k}, \ldots, \xi_{m,k}$ are associated for every $k$ and $(\xi_{1,k}, \ldots, \xi_{m,k}) \xrightarrow{\mathcal{D}} (\xi_1, \ldots, \xi_m)$ as $k \to \infty$, then $\xi_1, \ldots, \xi_m$ are associated.

(d) If two sets of associated variables are independent, then the union of these sets is also associated.

Then $\{T_{i,n}\}_{i=1}^{n-1}$ are associated for every $n$ by (a), (b) and (20). Analogously, $\{T_i^{(k)}\}_{i \in \mathbb{Z}}$ are associated for every $k$. Finally, since $T_i^{(k)} \to T_i$ a.s. as $k \to \infty$ for every $i$, (c) ensures the association of $\{T_i\}_{i \in \mathbb{Z}}$.

7. For any $s, t \in \mathbb{R}$ and $k \in \mathbb{Z}$,

$$\mathrm{cov}(\mathbb{1}_{\{T_0 \leq s\}}, \mathbb{1}_{\{T_k \leq t\}}) \geq 0. \tag{23}$$

*Proof.* This inequality follows from $\mathrm{cov}(\mathbb{1}_{\{T_0 \leq s\}}, \mathbb{1}_{\{T_k \leq t\}}) = \mathrm{cov}(\mathbb{1}_{\{s < T_0\}}, \mathbb{1}_{\{t < T_k\}})$, the association of $T_0, T_k$ and (b).

8. If $\mathbb{E} X_i^\gamma < \infty$ for some $\gamma \geq 2$, then the stationary sequence $\min\{T_i, t\}$ is strongly mixing for any $t < 1$ and its coefficients of strong mixing $\alpha(k)$ satisfy $\alpha(k) = o(k^{2-\gamma})$.

*Proof.* Recall that stationary r.v.'s $\xi_i$ are *strongly mixing* if $\alpha(k) \to 0$ as $k \to \infty$, where $\alpha(k)$ are the coefficients of strong mixing defined as

$$\alpha(k) := \sup_{A \in \mathcal{F}_{-\infty}^0, B \in \mathcal{F}_k^\infty} |\mathbb{P}(AB) - \mathbb{P}(A)\mathbb{P}(B)|;$$

here $\mathcal{F}_{-\infty}^0 := \sigma(\xi_0, \xi_{-1}, \ldots)$ and $\mathcal{F}_k^\infty := \sigma(\xi_k, \xi_{k+1}, \ldots)$ are the $\sigma$-algebras of "past" and "future," respectively. It is readily seen that

$$\alpha(k) \leq \sup_{0 \leq f, g \leq 1} |\mathrm{cov}(f(\xi_0, \xi_{-1}, \ldots), g(\xi_k, \xi_{k+1}, \ldots))|, \tag{24}$$

where the supremum is taken over Borel functions $f, g : \mathbb{R}^\infty \to [0, 1]$.

Let us estimate $\alpha(k)$ in the same way we estimated the left-hand side of (21). Fix some Borel functions $f, g : \mathbb{R}^\infty \to [0, 1]$. We approximate the variables from the "past" $T_0 \wedge t, T_{-1} \wedge t, T_{-2} \wedge t, \ldots$ by $T_0^{(k/2)} \wedge t, T_{-1}^{(k/2+1)} \wedge t, T_{-2}^{(k/2+2)} \wedge t, \ldots$, respectively; and for the variables from the "future," we use the analogous approximation. Now, $f(T_0^{(k/2)} \wedge t, T_{-1}^{(k/2+1)} \wedge t, \ldots)$ and $g(T_k^{(k/2)} \wedge t, T_{k+1}^{(k/2+1)} \wedge t, \ldots)$ are independent because the first is a function of $\{X_i\}_{i \leq k/2}$ and the second is a function of $\{X_i\}_{i \geq k/2+1}$. We then argue in the same way as in (22) to get

$$|\text{cov}(f(T_0 \wedge t, T_{-1} \wedge t, \ldots), g(T_k \wedge t, T_{k+1} \wedge t, \ldots))|$$
$$\leq 2\mathbb{P}\left\{\bigcup_{i=0}^{\infty}(T_{-i} \wedge t) \neq (T_{-i}^{(k/2+i)} \wedge t)\right\}$$
$$+ 2\mathbb{P}\left\{\bigcup_{i=0}^{\infty}(T_{k+i} \wedge t) \neq (T_{k+i}^{(k/2+i)} \wedge t)\right\}$$
$$\leq 4 \sum_{i=k/2}^{\infty} \mathbb{P}\{(T_0 \wedge t) \neq (T_0^{(i)} \wedge t)\}.$$

Now, by the formula of total probability, we have

$$\mathbb{P}\{(T_0 \wedge t) \neq (T_0^{(i)} \wedge t)\}$$
$$= \mathbb{P}\{(T_0 \wedge t) \neq (T_0^{(i)} \wedge t), T_0^{(i)} \geq t\} + \mathbb{P}\{(T_0 \wedge t) \neq (T_0^{(i)} \wedge t), T_0^{(i)} < t\}$$
$$\leq \mathbb{P}\{\mathbb{1}_{\{t \leq T_0\}} \neq \mathbb{1}_{\{t \leq T_0^{(i)}\}}\} + \mathbb{P}\{T_0 \neq T_0^{(i)}, T_0^{(i)} \leq t\}$$

and combining all the estimates together, by Lemma 1 (24) and arbitrariness of $f$ and $g$, we get $\alpha(k) \leq 8 \sum_{i=k/2}^{\infty} o(i^{1-\gamma}) = o(k^{2-\gamma})$ if $\gamma > 2$. For $\gamma = 2$, we get $\alpha(k) \leq 16 \sum_{i=k/2}^{\infty} \mathbb{P}\{\inf_{i \geq M} \frac{S_i}{i} < t^2\} = o(1)$ using the same argument and applying (13), (14), and Fact 2 instead of Lemma 1.

3.5. *The last collision.* We finish this section with a statement on the convergence of the moments of the last collision.

PROPOSITION 3. *In the i.d. model,* $T_n^{\text{last}} \xrightarrow{\mathbb{P}} 1$ *as* $n \to \infty$ *if* $\mathbb{E}X_i^2 < \infty$.

This result is well known for the Poisson model; see Giraud [8].

PROOF OF PROPOSITION 3. Let us first prove that $\mathbb{P}\{T_n^{\text{last}} \geq t\} \to 0$ as $n \to \infty$ for all $t > 1$. Since $T_n^{\text{last}} = \max_{1 \leq j \leq n-1} T_{j,n}$, we have

$$(25) \qquad \mathbb{P}\{T_n^{\text{last}} \geq t\} \leq \sum_{j=1}^{n-1} \mathbb{P}\{T_{j,n} \geq t\}.$$

By taking into account that the minima in (17) are nonpositive and by arguing as in (18),

$$\mathbb{P}\{T_{j,n} \geq t\} \leq \mathbb{P}\left\{\min_{1 \leq k \leq j \vee n-j} \frac{1}{k}\sum_{i=1}^{k-1}(k-i)(X_{j+i+1}-t^2) + X_{j+1} - t^2 \geq 0\right\}$$

$$= \mathbb{P}\left\{\min_{1 \leq k \leq j \vee n-j} \sum_{i=1}^{k}(k-i+1)(X_i - t^2) \geq 0\right\}$$

$$\leq \mathbb{P}\left\{\min_{1 \leq k \leq n/2} \sum_{i=1}^{k}(S_i - it^2) \geq 0\right\}.$$

We claim that (without any assumptions on the moments of $X_i$)

(26) $$\mathbb{P}\{T_{j,n} \geq t\} \leq \mathbb{P}\left\{\sup_{i \geq (t-1)/4tn} \frac{S_i}{i} > \frac{1+t^2}{2}\right\};$$

recall that $t > 1$. Clearly, (26) follows if we check that

$$\left\{\min_{1 \leq k \leq n/2}\sum_{i=1}^{k}(S_i - it^2) \geq 0\right\} \subset \left\{\sup_{i \geq (t-1)/4tn} \frac{S_i}{i} > \frac{1+t^2}{2}\right\}.$$

Assume the converse; then, by the nonnegativity of $S_i$,

$$0 \leq \sum_{i=1}^{n/2}(S_i - it^2) = \sum_{i=1}^{cn}(S_i - it^2) + \sum_{i=cn+1}^{n/2}(S_i - it^2)$$

$$\leq \sum_{i=1}^{cn}(S_{cn} - it^2) + \sum_{i=cn+1}^{n/2}\left(i\frac{1+t^2}{2} - it^2\right),$$

where $c := \frac{t-1}{4t}$. We estimate the last expression with

$$cnS_{cn} - \frac{(cn)^2}{2}t^2 - \frac{(n/2)^2 - (cn)^2}{2} \cdot \frac{t^2-1}{2} \leq \frac{c^2}{2}n^2 - \frac{1/4 - c^2}{2} \cdot \frac{t^2-1}{2}n^2.$$

It is simple to check that the right-hand side is negative, thus we have a contradiction.

Then from (25), (26) and Fact 2 it follows that $\mathbb{P}\{T_n^{\text{last}} \geq t\} = \sum_{i=1}^{n-1} o((cn)^{-1}) = o(1)$ for all $t > 1$.

Now let us prove that $\mathbb{P}\{T_n^{\text{last}} < t\} \to 0$ as $n \to \infty$ for all $t < 1$. Since $T_n^{\text{last}} = \max_{1 \leq j \leq n-1} T_{j,n}$, we estimate

$$\mathbb{P}\{T_n^{\text{last}} < t\} \leq \mathbb{P}\left\{\max_{1 \leq j \leq \sqrt{n}-1} T_{j\sqrt{n},n} < t\right\}$$

$$= \mathbb{P}\left\{\max_{1 \leq j \leq \sqrt{n}-1} T_{j\sqrt{n}}^{(\sqrt{n}/2)} < t\right\} + \sum_{j=1}^{\sqrt{n}-1}\mathbb{P}\{\mathbb{1}_{\{t \leq T_{j\sqrt{n},n}\}} \neq \mathbb{1}_{\{t \leq T_{j\sqrt{n}}^{(\sqrt{n}/2)}\}}\}.$$

In view of (10) and Lemma 1, the sum is $\sum_{j=1}^{\sqrt{n}-1} o(n^{-1/2}) = o(1)$, hence it remains to check that the first probability in the last line tends to zero. For a fixed $n$, all $T_{j\sqrt{n}}^{(\sqrt{n}/2)}$ are independent because each one is a function of $\{X_i\}_{|j\sqrt{n}-i|\leq\sqrt{n}/2}$ (to be precise, of $X_{j\sqrt{n}-\sqrt{n}/2+2},\ldots,X_{j\sqrt{n}+\sqrt{n}/2}$). Thus,

$$\mathbb{P}\left\{\max_{1\leq j\leq\sqrt{n}-1} T_{j\sqrt{n}}^{(\sqrt{n}/2)} < t\right\} = \mathbb{P}^{\sqrt{n}-1}\{T_{\sqrt{n}}^{(\sqrt{n}/2)} < t\} \leq \mathbb{P}^{\sqrt{n}-1}\{T_0 < t\},$$

which tends to zero; indeed, $\mathbb{P}\{T_0 < t\} < 1$ by Property 3, Section 3.4. □

**4. Proofs of Fact 1 and Theorem 1 for the i.d. model.** Recall that the number of clusters $K_n(t)$ is given by (6). Our idea is to study $\sum_{i=1}^{n-1} \mathbb{1}_{\{t<T_i\}}$ instead of $\sum_{i=1}^{n-1} \mathbb{1}_{\{t<T_{i,n}\}}$: We thus deal with a single sequence $T_i$ and avoid considering the triangular array $T_{i,n}$.

Let us now prove Fact 1 for the i.d. model. We prove (1) for $t \neq 1$ without any additional assumptions on $X_i$; for $t = 1$, we require $\mathbb{E}X_i^2 < \infty$. The properties of the limit function $a(t)$ were studied in Section 3.4, Properties 3 and 4.

PROOF OF FACT 1. We put
$$(27) \qquad a(t) := \mathbb{P}\{T_0 > t\}.$$

Let us first prove (1) for all $t < 1$. It is sufficient to check that

$$(28) \qquad \frac{K_n(t)}{n} - \frac{1}{n}\sum_{i=1}^n \mathbb{1}_{\{t<T_i\}} \xrightarrow{\mathbb{P}} 0, \qquad n \to \infty.$$

Indeed, the stationary sequence $\mathbb{1}_{\{t<T_i\}}$ satisfies the law of large numbers by Property 5, Section 3.4, and the well-known result of S. N. Bernstein:

FACT 3. *The law of large numbers holds for r.v.'s $\xi_i$ if there exists a sequence $r(k) \to 0$ such that $\operatorname{cov}(\xi_i, \xi_j) \leq r(|i-j|)$ for all $i, j \in \mathbb{N}$.*

By (6),
$$\left|\frac{K_n(t)}{n} - \frac{1}{n}\sum_{i=1}^n \mathbb{1}_{\{t<T_i\}}\right| \leq \frac{1}{n} + \frac{1}{n}\sum_{i=1}^{n-1}(\mathbb{1}_{\{t<T_{i,n}\}} - \mathbb{1}_{\{t<T_i\}}),$$

where we used (10) to get the nonnegativity of the right-hand side. Then (28) immediately follows from the Chebyshev inequality provided that the expectation of the right-hand side tends to zero. By using (10), we obtain

$$\frac{1}{n}\sum_{i=1}^{n-1}\mathbb{E}(\mathbb{1}_{\{t<T_{i,n}\}} - \mathbb{1}_{\{t<T_i\}}) \leq \frac{1}{n}\sum_{i=1}^{n-1}(\mathbb{E}\mathbb{1}_{\{t<T_i^{(i\wedge n-i)}\}} - \mathbb{E}\mathbb{1}_{\{t<T_i\}})$$
$$= \frac{1}{n}\sum_{i=1}^{n-1}\mathbb{P}\{\mathbb{1}_{\{t<T_i\}} \neq \mathbb{1}_{\{t<T_i^{(i\wedge n-i)}\}}\},$$

which is $\frac{2}{n}\sum_{i=1}^{n/2} o(1) = o(1)$ by Lemma 1. To be very precise, Lemma 1 deals with slightly different indicators, but we can estimate the considered probability by repeating the proof of Lemma 1 word for word (or just use Property 4, Section 3.4).

We now check that (1) holds for all $t > 1$. Using (26) gives $\mathbb{E}\frac{K_n(t)}{n} = \frac{1}{n}\sum_{i=1}^{n-1} \mathbb{P}\{T_{i,n} > t\} \to 0$ as $n \to \infty$ and $\frac{K_n(t)}{n} \xrightarrow{\mathbb{P}} a(t) = 0$ follows from the Chebyshev inequality.

It remains to check that (1) holds for $t = 1$ if $\mathbb{E}X_i^2 < \infty$ to conclude the proof. If $\mathbb{D}X_i = 0$, then the situation is deterministic, this case was described in Introduction. Here we always have $K_n(1) = 1$ and (1) is true. If $0 < \mathbb{D}X_i < \infty$, then by Property 3 from Section 3.4, we have $a(1) = 0$ and $\mathbb{P}\{T_0 = 1\} = 0$; consequently, $a(t) = \mathbb{P}\{T_0 > t\}$ is continuous at $t = 1$. Then (1) is true for $t = 1$ since $0 < \frac{K_n(1)}{n} \le \frac{K_n(t)}{n} \xrightarrow{\mathbb{P}} a(t)$ for any $t \in (0,1)$ and $a(t) \to a(1) = 0$ as $t \nearrow 1$. $\square$

Now we prove Theorem 1 for the i.d. model. We think of $D[0,1]$ as of a separable metric space equipped with the Skorohod metric $d$, which induces the Skorohod topology.

PROOF OF THEOREM 1. At first, we prove (2). In view of representation (6) for $K_n(t)$, relation (2) follows from the relation

$$(29) \quad \sup_{0 \le t \le 1-\varepsilon} \left| \frac{1}{\sqrt{n}} \sum_{i=1}^{n-1} \mathbb{1}_{\{t < T_{i,n}\}} - \frac{1}{\sqrt{n}} \sum_{i=1}^{n} \mathbb{1}_{\{t < T_i\}} \right| \xrightarrow{\mathbb{P}} 0 \quad \text{for all } \varepsilon \in (0,1)$$

and the existence of a centered Gaussian process $K(\cdot)$ on $[0,1)$ such that

$$(30) \quad \frac{1}{\sqrt{n}}\left\{ \sum_{i=1}^{n} \mathbb{1}_{\{t < T_i\}} - na(t) \right\} \xrightarrow{\mathcal{D}} K(\cdot) \quad \text{in } D[0, 1-\varepsilon] \text{ for all } \varepsilon \in (0,1).$$

Indeed, if $Y_n \xrightarrow{\mathcal{D}} Y$ and $d(Y_n, Y_n') \xrightarrow{\mathbb{P}} 0$ for some random elements $Y_n, Y_n', Y$ of the separable metric space $D[0, 1-\varepsilon]$, then $Y_n' \xrightarrow{\mathcal{D}} Y$; recall that $d(Y_n, Y_n') \le \sup_{t \in [0, 1-\varepsilon]} |Y_n(t) - Y_n'(t)|$.

We start with (29). It is sufficient to prove that the expectation of the left-hand side tends to zero. Since the supremum of a sum does not exceed the sum of suprema, let us check that

$$(31) \quad \frac{1}{\sqrt{n}} \sum_{i=1}^{n-1} \mathbb{E} \sup_{0 \le t \le 1-\varepsilon} |\mathbb{1}_{\{t<T_{i,n}\}} - \mathbb{1}_{\{t<T_i\}}| \longrightarrow 0 \quad \text{for all } \varepsilon \in (0,1).$$

By (10), we have

$$\mathbb{E} \sup_{0 \le t \le 1-\varepsilon} |\mathbb{1}_{\{t<T_{i,n}\}} - \mathbb{1}_{\{t<T_i\}}| \le \mathbb{E} \sup_{0 \le t \le 1-\varepsilon} (\mathbb{1}_{\{t<T_i^{(i \wedge n-i)}\}} - \mathbb{1}_{\{t<T_i\}})$$

$$= \mathbb{P}\{T_i \neq T_i^{(i \wedge n-i)}, T_i \leq 1-\varepsilon\}$$
$$= \mathbb{P}\{T_i \neq T_i^{(i \wedge n-i)}, T_i^{(i \wedge n-i)} < 1-\varepsilon\}$$
$$+ \mathbb{P}\{\mathbb{1}_{\{1-\varepsilon \leq T_i\}} \neq \mathbb{1}_{\{1-\varepsilon \leq T_i^{(i \wedge n-i)}\}}\},$$

where the last equality was obtained via the formula of total probability. Combining the estimates together and using Lemma 1,

$$\frac{1}{\sqrt{n}} \sum_{i=1}^{n-1} \mathbb{E} \sup_{0 \leq t \leq 1-\varepsilon} |\mathbb{1}_{\{t<T_{i,n}\}} - \mathbb{1}_{\{t<T_i\}}|$$
$$\leq \frac{2\rho(1-\varepsilon)}{\sqrt{n}} \sum_{i=1}^{n-1} (i \wedge n-i)^{1-\gamma} = \frac{4\rho(1-\varepsilon)}{\sqrt{n}} \sum_{i=1}^{n/2} i^{1-\gamma}.$$

The last expression is $O(n^{3/2-\gamma})$ and (31), which implies (29), follows.

Now let us prove (30). As long as

$$U_n(t) := -\frac{1}{\sqrt{n}} \left\{ \sum_{i=1}^n \mathbb{1}_{\{t<T_i\}} - na(t) \right\} = \sqrt{n} \left\{ \frac{1}{n} \sum_{i=1}^n \mathbb{1}_{\{T_i \leq t\}} - (1-a(t)) \right\},$$

the $U_n(\cdot)$ is the empirical process of stationary r.v.'s $T_i$ with the continuous common distribution function $1-a(t)$. By $K(\cdot) \stackrel{\mathcal{D}}{=} -K(\cdot)$, (30) is equivalent to the existence of a centered Gaussian process $K(\cdot)$ on $[0,1)$ such that

(32) $\qquad U_n(\cdot) \xrightarrow{\mathcal{D}} K(\cdot) \qquad$ in $D[0, 1-\varepsilon]$ for all $\varepsilon \in (0,1)$.

We will use the following result from Lin and Lu [17], Section 12 on convergence of empirical processes. They attribute this statement to Q.-M. Shao, who published it in 1986, in Chinese.

FACT 4. *Let $\xi_i$ be a sequence of stationary strongly mixing r.v.'s distributed on $[0,1]$, and let $F$ be the common distribution function of $\xi_i$. Suppose $F(x) = x$ on $[0,1]$ (i.e., $\xi_i$ are uniformly distributed) and the coefficients of strong mixing of the sequence $F(\xi_i)$ decrease as $O(k^{-(2+\delta)})$ as $k \to \infty$ for some $\delta > 0$. Then the empirical processes of $\xi_i$ weakly converge in $D[0,1]$ to a centered Gaussian process with the covariance function $\sum_{i \in \mathbb{Z}} \mathrm{cov}(\mathbb{1}_{\{\xi_0 \leq s\}}, \mathbb{1}_{\{\xi_i \leq t\}})$.*

REMARK. The limit Gaussian process is a.s. continuous on $[0,1]$. Fact 4 also holds true if $F$ is an arbitrary continuous distribution function.

The a.s. continuity of the limit process could be concluded by a comparison of the proof from Lin and Lu [17] with the proof of Theorem 22.1 from Billingsley [3]. The statements and the proofs of these theorems are identical, but Lin and Lu do not state the continuity while Billingsley does. Further, since $F(\xi_i)$ is uniformly distributed on $[0,1]$ if $F$ is continuous, Fact 4 holds true for every continuous $F$; see the proof of Theorem 22.1 by Billingsley [3] for explanations.

Recall that we need to prove the convergence of the empirical process of $T_i$. It seems that the r.v.'s $T_i$ are not strongly mixing; but $\min\{T_i, 1-\varepsilon\}$ are strongly mixing because of Property 8, Section 3.4. These variables are not continuous and so we need to fix them. Let us fix an $\varepsilon \in (0,1)$, and let $\alpha_i$ be i.i.d. r.v.'s independent of all $T_i$ and, say, uniformly distributed on $[0,\varepsilon]$; we define $\tilde{T}_i := \min\{T_i, 1-\varepsilon\} + \mathbb{1}_{\{T_i \geq 1-\varepsilon\}} \alpha_i$.

The stationary variables $\tilde{T}_i$ are distributed on $[0,1]$, their common distribution function $G$ is continuous, and the coefficients of strong mixing of $G(\tilde{T}_i)$ decrease as $o(k^{2-\gamma})$. The proof of the last statement is the same as the proof of Property 8 from Section 3.4. Indeed, approximate the variables $G(\tilde{T}_0), G(\tilde{T}_{-1}), \ldots$ from the "past" by $G(\tilde{T}_0^{(k/2)}), G(\tilde{T}_{-1}^{(k/2+1)}), \ldots$ where $\tilde{T}_i^{(m)} := \min\{T_i^{(m)}, 1-\varepsilon\} + \mathbb{1}_{\{T_i^{(m)} \geq 1-\varepsilon\}} \alpha_i$; use the analogous approximation for the variables from the "future"; and then repeat word for word the arguments of the previous proof.

Now, recalling that $\gamma > 4$, we see that $\tilde{T}_i$ satisfy the assumptions of Fact 4, with the only difference that their distribution is not uniform. By $\tilde{U}_n(\cdot)$ denote the empirical process of $\tilde{T}_i$; clearly, $\tilde{U}_n(\cdot)$ coincides with the empirical process $U_n(\cdot)$ of $T_i$ on $[0, 1-\varepsilon]$. By the remark to Fact 4, we conclude that first,

$$\tilde{U}_n(\cdot) \xrightarrow{\mathcal{D}} \tilde{K}(\cdot) \qquad \text{in } D[0,1], \tag{33}$$

where $\tilde{K}(\cdot)$ is a centered Gaussian process with the covariance function

$$\tilde{R}(s,t) := \sum_{i \in \mathbb{Z}} \text{cov}(\mathbb{1}_{\{\tilde{T}_0 \leq s\}}, \mathbb{1}_{\{\tilde{T}_i \leq t\}})$$

and, second, trajectories of $\tilde{K}(\cdot)$ are a.s. continuous on $[0,1]$.

[There exists a simpler and more elegant proof of (33). Note that $\{\tilde{T}_i\}_{i \in \mathbb{Z}}$ are associated as coordinate-wise nondecreasing functions of associated r.v.'s $\{T_i, \alpha_i\}_{i \in \mathbb{Z}}$, see (a), (b) and (d) from Property 6, Section 3.4. Then we can obtain (33) applying the result of Louhichi [18] on convergence of empirical processes of stationary associated r.v.'s $\xi_i$ instead of using Fact 4. This theorem requires only $\text{cov}(F(\xi_0), F(\xi_k)) = O(k^{-(4+\delta)})$, which could be proved analogously to Property 5, Section 3.4. Thus we avoid the complicated estimations of the strong mixing coefficients, and the proof of (33) is becomes much simpler. The only problem is that this proof requires $\gamma > 5$.

We also note that the a.s. continuity of $\tilde{K}(\cdot)$ could be proved directly, without referring to the proof of Fact 4. The arguments should be the same as in the proof of the continuity of $K^{\text{Unif}}(\cdot)$ in Section 5.]

Define

$$R(s,t) := \sum_{i \in \mathbb{Z}} \text{cov}(\mathbb{1}_{\{T_0 \leq s\}}, \mathbb{1}_{\{T_i \leq t\}}), \tag{34}$$

which is, evidently, equal to $\tilde{R}(s,t)$ on $[0, 1-\varepsilon]^2$. Since $\tilde{R}(s,t)$ is positive definite and $\varepsilon > 0$ is arbitrary, the function $R(s,t)$ is positive definite on $[0,1)^2$. Therefore, by Lifshits [15], Section 4, there exists a centered Gaussian process $K(\cdot)$ on $[0,1)$ with the covariance function $R(s,t)$. The trajectories of $K(\cdot)$ are a.s. continuous on $[0,1)$ by $K(\cdot) \stackrel{\mathcal{D}}{=} \tilde{K}(\cdot)$ on $[0, 1-\varepsilon]$, arbitrariness of $\varepsilon > 0$, and the a.s. continuity of $\tilde{K}(\cdot)$ on $[0,1]$.

Finally, by (33), $\tilde{U}_n(\cdot) = U_n(\cdot)$ on $[0, 1-\varepsilon]$, $\tilde{K}(\cdot) \stackrel{\mathrm{D}}{=} K(\cdot)$ on $[0, 1-\varepsilon]$, and the a.s. continuity of $\tilde{K}(\cdot)$, we get (32). Since (32) implies (30), we conclude the proof of (2).

Only the stated properties of $R(s,t)$ remain to be proven. The continuity of the joint distribution function of continuous variables $T_0$ and $T_i$ implies that $\text{cov}(\mathbb{1}_{\{T_0 \leq s\}}, \mathbb{1}_{\{T_i \leq t\}})$ is continuous on $[0,1)^2$ for every $i \geq 0$. Then, in view of (21), $R(s,t)$ is continuous on $[0,1)^2$ as a sum of uniformly converging series of continuous functions.

The strict positivity of $R(s,t)$ on $(\sqrt{\mu}, 1)^2$ trivially follows from (34), (23) and $\text{cov}(\mathbb{1}_{\{T_0 \leq s\}}, \mathbb{1}_{\{T_0 \leq t\}}) = a(s \vee t)(1 - a(s \wedge t)) > 0$; the last inequality holds by Property 3, Section 3.4. The $R(s,t) = 0$ on $[0,1)^2 \setminus (\sqrt{\mu}, 1)^2$ follows from $\mathbb{P}\{T_i \leq \sqrt{\mu}\} = 0$, which holds by Properties 3 and 4 from Section 3.4. □

We note that (3) holds for $t \neq 1$ under the less restrictive condition $\mathbb{E}X_i^2 < \infty$. For $t < 1$, the proof is almost the same: By (29), which is true for $\gamma > 3/2$, we conclude that (3) holds if the stationary associated sequence $\mathbb{1}_{\{t < T_i\}}$ satisfies the central limit theorem. Then we refer to the central limit theorem for stationary associated sequences from Newman [21]; his theorem requires only $R(t,t) < \infty$, that is, the convergence of the right-hand side of (34). This condition holds by (13) and Fact 2. For $t > 1$, relation (3) holds true with $\sigma^2(t) = 0$ because of Proposition 3.

Finally, note that the process $K(\cdot)$ is associated, that is, the r.v.'s $\{K(t)\}_{t \in [0,1)}$ are associated. In fact, by (6), Property 6 from Section 3.4, and Condition (b) from the same Property 6, the processes $\frac{K_n(\cdot) - na(\cdot)}{\sqrt{n}}$ are associated for every $n$. Then $K(\cdot)$ is associated by (2) and (c), Property 6.

**5. Proof of Theorem 1 for the uniform model.** There exists a simple method that allows to extend results from the Poisson model to the uniform model and vise versa. The method is based on the next statement (see Karlin [13], Section 9.1).

FACT 5. *Let $S_i$ be an exponential random walk. Then for any $k \geq 1$, we have*

$$(35) \qquad \left(\frac{S_1}{S_{k+1}}, \frac{S_2}{S_{k+1}}, \ldots, \frac{S_k}{S_{k+1}}\right) \stackrel{\mathcal{D}}{=} (U_{1,k}, U_{2,k}, \ldots, U_{k,k}),$$

*where $U_{i,k}$ are the order statistics of $k$ i.i.d. r.v.'s uniformly distributed on $[0,1]$. Moreover, the random vector in the left-hand side of (35) is independent of $S_{k+1}$.*

Therefore, if $x_{j,n}^{\text{Poiss}}(0) = \frac{1}{n} S_j$ are the initial positions of particles in the Poisson model, then for the initial positions of particles in the uniform model, we have $x_{j,n}^{\text{Unif}}(0) = \frac{n}{S_{n+1}} \cdot x_{j,n}^{\text{Poiss}}(0)$. By Proposition 2 and (5), we conclude that

$$(36) \qquad T_{j,n}^{\text{Unif}} = \beta_n^{-1} T_{j,n}^{\text{Poiss}}, \qquad \beta_n := \sqrt{\frac{S_{n+1}}{n}},$$

and hence, using (6), we get

$$(37) \qquad K_n^{\text{Unif}}(t) = K_n^{\text{Poiss}}(\beta_n t).$$

Note that the process $K_n^{\text{Unif}}(\cdot)$ and the r.v. $\beta_n$ are independent since values of the process are defined by $x_{1,n}^{\text{Unif}}(0), \ldots, x_{n,n}^{\text{Unif}}(0)$, which are mutually independent of $\beta_n$ by Fact 5.

Now we prove Theorem 1 for the uniform model.

PROOF OF THEOREM 1. Denote

$$Y_n(t) := \frac{K_n^{\text{Unif}}(t) - na(t)}{\sqrt{n}}, \qquad Z_n(t) := \sqrt{n}(a(t) - a(\beta_n t));$$

we stress that $Y_n(\cdot)$ and $Z_n(\cdot)$ are independent.

Fix an $\varepsilon \in (0, 1)$. First, it follows from (2) for the Poisson model and (37) that

$$(38) \qquad Y_n(\cdot) + Z_n(\cdot) \xrightarrow{\mathcal{D}} K^{\text{Poiss}}(\cdot) \qquad \text{in } D[0, 1-\varepsilon].$$

Indeed, the process $Y_n(\cdot) + Z_n(\cdot)$ is obtained from $\frac{1}{\sqrt{n}}(K_n^{\text{Poiss}}(\cdot) - na(\cdot))$ by the random time change $t \mapsto \beta_n t$; and since $\|\beta_n t - t\|_{C[0,1-\varepsilon]} \xrightarrow{\mathbb{P}} 0$, we have

$$d\left(Y_n(\cdot) + Z_n(\cdot), \frac{K_n^{\text{Poiss}}(\cdot) - na(\cdot)}{\sqrt{n}}\right) \xrightarrow{\mathbb{P}} 0$$

by the definition of the Skorohod metric $d$.

Second, from Fact 1, (15), and (27) it follows that $a^{\text{Unif}}(t) = a^{\text{Poiss}}(t) = \mathbb{P}\{T_0^{\text{Poiss}} \geq t\} = 1 - t^2$ for $0 \leq t \leq 1$, and by the central limit theorem,

$$(39) \qquad Z_n(t) \xrightarrow{\mathcal{D}} t^2 \eta \qquad \text{in } D[0, 1-\varepsilon],$$

where $\eta$ is a standard Gaussian r.v.

We claim that (38), the independence of $Y_n(\cdot)$ and $Z_n(\cdot)$, and (39) yield the weak convergence of $Y_n(\cdot)$ in $D[0, 1-\varepsilon]$. Let us check the tightness of $Y_n(\cdot)$ and the convergence of their finite-dimensional distributions.

The tightness of $Y_n(\cdot)$ in $D[0, 1-\varepsilon]$ follows from $Y_n(\cdot) = (Y_n(\cdot) + Z_n(\cdot)) - Z_n(\cdot)$, (38), and (39). Indeed, by the Prokhorov theorem, (38) and (39) yield that both sequences $Y_n(\cdot) + Z_n(\cdot)$ and $-Z_n(\cdot)$ are tight. But trajectories of $-Z_n(\cdot)$ are a.s. continuous because of the continuity of $a(\cdot)$, and the tightness follows from the continuity of addition $+ : D \times C \to D$ and the fact that under any continuous mapping, the image of a compact set is also a compact set.

Now we study convergence of finite dimensional distributions of $Y_n(\cdot)$. Recall that the characteristic function of a centered Gaussian vector in $\mathbb{R}^m$ is $e^{-1/2(R\mathbf{u}, \mathbf{u})}$, where $\mathbf{u} \in \mathbb{R}^m$ and $R$ is the covariance matrix of the vector. Then (38), the independence of $Y_n(\cdot)$ and $Z_n(\cdot)$, and (39) yield that for the characteristic functions of *all* finite-dimensional distributions of $Y_n(\cdot)$, we have

$$(40) \qquad \mathbb{E} e^{i(Y_n(\mathbf{t}), \mathbf{u})} \longrightarrow e^{-1/2(\{R^{\mathrm{Poiss}}(t_j, t_k) - t_j^2 t_k^2\}_{j,k=1}^m \mathbf{u}, \mathbf{u})},$$

where $\mathbf{u} \in \mathbb{R}^m, \mathbf{t} = (t_1, \ldots, t_m) \in [0, 1-\varepsilon]^m$, and $Y_n(\mathbf{t}) := (Y_n(t_1), \ldots, Y_n(t_m))$. We stress that (40) is true for every $\mathbf{t} \in [0, 1-\varepsilon]^m$ since the limit processes in (38) and (39) have continuous trajectories.

We see that the matrix $\{R^{\mathrm{Poiss}}(t_j, t_k) - t_j^2 t_k^2\}_{j,k=1}^m$ is positive definite for any $\mathbf{t} = (t_1, \ldots, t_m) \in [0, 1-\varepsilon]^m$ and $m \geq 1$ since the absolute value of the left-hand side of (40) does not exceed one. Putting

$$R^{\mathrm{Unif}}(s, t) := R^{\mathrm{Poiss}}(s, t) - s^2 t^2,$$

we have $\{R^{\mathrm{Poiss}}(t_j, t_k) - t_j^2 t_k^2\}_{j,k=1}^m = \{R^{\mathrm{Unif}}(t_j, t_k)\}_{j,k=1}^m$; then the function $R^{\mathrm{Unif}}(s, t)$ is positive definite on $[0, 1)^2$ since $\varepsilon > 0$ is arbitrary. Thus, by Lifshits [15], Section 4, $R^{\mathrm{Unif}}(s, t)$ is the covariance function of some centered Gaussian process $K^{\mathrm{Unif}}(\cdot)$ on $[0, 1)$.

Relation (2) is thus proved. Now check that $K^{\mathrm{Unif}}(\cdot) \in C[0, 1-\varepsilon]$ a.s. to conclude the proof of Theorem 1 for the uniform model.

For this purpose, let us prove that a.s., trajectories of $Y_n(\cdot)$ have jumps of size $\frac{1}{\sqrt{n}}$ only. In fact, the jumps of $Y_n(\cdot)$ coincide with the jumps of $\frac{1}{\sqrt{n}} K_n^{\mathrm{Unif}}(\cdot)$, whose jumps are of size $\frac{1}{\sqrt{n}}$ if and only if $T_{j_1,n}^{\mathrm{Unif}} \neq T_{j_2,n}^{\mathrm{Unif}}$ for $1 \leq j_1 \neq j_2 \leq n-1$. By (36), we need to verify that $T_{j_1,n}^{\mathrm{Poiss}} \neq T_{j_2,n}^{\mathrm{Poiss}}$ a.s. for $1 \leq j_1 \neq j_2 \leq n-1$. This relation follows from (20) if $H(k_1, j_1, l_1) \neq H(k_2, j_2, l_2)$ a.s. for $j_1 \neq j_2$ and $k_1, k_2, l_1, l_2 \geq 1$. The last a.s. nonequality is obvious because if the equality holds true, then a certain nontrivial linear combination of i.i.d. exponential $X_i$ equals zero.

Then there exist a.s. continuous $\tilde{Y}_n(\cdot)$ such that $\sup_{t\in[0,1-\varepsilon]}|\tilde{Y}_n(t)-Y_n(t)| \leq \frac{1}{\sqrt{n}}$ a.s.; consequently, $d(\tilde{Y}_n, Y_n) \leq \frac{1}{\sqrt{n}}$ a.s. Then by $Y_n(\cdot) \xrightarrow{\mathcal{D}} K^{\mathrm{Unif}}(\cdot)$, we also have $\tilde{Y}_n(\cdot) \xrightarrow{\mathcal{D}} K^{\mathrm{Unif}}(\cdot)$. But $1 = \liminf \mathbb{P}\{\tilde{Y}_n(\cdot) \in C\} \leq \mathbb{P}\{K^{\mathrm{Unif}}(\cdot) \in C\}$ since $C \subset D$ is closed in the Skorohod topology, therefore, a.s., $K^{\mathrm{Unif}}(\cdot)$ is continuous on $[0, 1-\varepsilon]$.

Since $\varepsilon \in (0,1)$ is arbitrary, a.s., $K^{\mathrm{Unif}}(\cdot)$ is continuous on the whole interval $[0,1)$. The $R^{\mathrm{Unif}}(s,t) = R^{\mathrm{Poiss}}(s,t) - s^2 t^2$ is continuous on $[0,1)^2$ because $R^{\mathrm{Poiss}}(s,t)$ is. $\square$

**6. The number of clusters at the critical moment.** Now we turn our attention to the number of clusters at the critical moment $t = 1$. We are interested in the behavior of

$$\frac{K_n(1) - na(1)}{\sqrt{n}} = \frac{K_n(1)}{\sqrt{n}},$$

which is the left-hand side of (3) at $t = 1$; here we have $a(1) = 0$ under $\mathbb{E}X_i^2 < \infty$, see Property 3, Section 3.4.

We do not know if this sequence is weakly convergent, but we hope that it is. We also have a naive guess that its limit is Gaussian because the limit in Theorem 1 is Gaussian. In view of $K_n(1) \geq 1$, this conjectured weak limit is nonnegative, hence it is Gaussian if and only if it is identically equal to zero. However, the results of this section show that the limit is nonzero, thus our guess on Gaussianity fails.

The study of convergence of $\frac{K_n(1)}{\sqrt{n}}$ is quite complicated. Therefore, in this section, we consider only the Poisson model. First, let us prove the following statement.

PROPOSITION 4. *In the Poisson model, we have* $\lim_{n\to\infty} \mathbb{P}\{K_n(1) = 1\} > 0$.

PROOF. On the one hand, $K_n(1) = 1$ is equivalent to $T^{\mathrm{last}}_{n;\mathrm{Poiss}} \leq 1$, where $T^{\mathrm{last}}_{n;\mathrm{Poiss}}$ denotes the moment of the last collision in the Poisson model. On the other hand, a result by Giraud [8] states that in the uniform model,

$$\sqrt{n}(T^{\mathrm{last}}_{n;\mathrm{Unif}} - 1) \xrightarrow{\mathcal{D}} \sup_{0\leq x \leq 1}\left(\frac{1}{1-x}\int_x^1 \overset{\circ}{W}(y)\,dy - \frac{1}{x}\int_0^x \overset{\circ}{W}(y)\,dy\right) =: \tau,$$

where $\overset{\circ}{W}(\cdot)$ is a Brownian bridge. Now, by (36), we have $T^{\mathrm{last}}_{n;\mathrm{Unif}} = \beta_n^{-1} T^{\mathrm{last}}_{n;\mathrm{Poiss}}$, hence

(41) $$\sqrt{n}(\beta_n^{-1} T^{\mathrm{last}}_{n;\mathrm{Poiss}} - 1) \xrightarrow{\mathcal{D}} \tau.$$

But from the central limit theorem and the law of large numbers,

$$(42) \quad \sqrt{n}(\beta_n^{-1} - 1) = -\frac{S_{n+1} - n}{\sqrt{n}} \cdot \frac{n}{\sqrt{S_{n+1}}(\sqrt{S_{n+1}} + \sqrt{n})} \xrightarrow{\mathcal{D}} \frac{\eta}{2},$$

where $\eta$ is a standard Gaussian r.v. and $S_i$ is a standard exponential random walk that defines initial positions of particles. Since, in view of Fact 5, $T_{n;\mathrm{Unif}}^{\mathrm{last}} = \beta_n^{-1} T_{n;\mathrm{Poiss}}^{\mathrm{last}}$ and $\beta_n$ are independent, from (41), (42), and the law of large numbers it follows that

$$\sqrt{n}(T_{n;\mathrm{Poiss}}^{\mathrm{last}} - 1) \xrightarrow{\mathcal{D}} \tau - \frac{\eta}{2} \stackrel{\mathcal{D}}{=} \tau + \frac{\eta}{2},$$

where $\tau$ and $\eta$ are independent. Thus,

$$\lim_{n \to \infty} \mathbb{P}\{K_n(1) = 1\} = \lim_{n \to \infty} \mathbb{P}\{T_{n;\mathrm{Poiss}}^{\mathrm{last}} \leq 1\} = \mathbb{P}\left\{\tau + \frac{\eta}{2} \leq 0\right\} > 0. \quad \square$$

The main advantage of the Poisson model is that, by Lemma 2 and Property 4, Section 3.4 we have $\mathbb{P}\{T_{j,n} > 1\} = ep_j p_{n-j}$, where

$$p_k := \mathbb{P}\left\{\min_{1 \leq m \leq k} \sum_{i=1}^{m} (S_i - \mathbb{E}S_i) \geq 0\right\}$$

and $S_i$ is a standard exponential random walk. We say that the sequence of r.v.'s $\sum_{i=1}^{m}(S_i - \mathbb{E}S_i)$ is an *integrated random walk*. In the proof of Property 3, Section 3.4, we showed that $p_k \to 0$ as $k \to \infty$. Therefore, it is reasonable to say that $p_k$ are the *unilateral small deviation probabilities of an integrated centered random walk*.

We need to obtain the asymptotics of $p_k \to 0$ to continue the study of convergence of $\frac{K_n(1)}{\sqrt{n}}$. Unfortunately, the results of the rest of this section are completely dependent on the correctness of the following conjecture.

CONJECTURE 1. *We have $p_k \sim c_1 k^{-1/4}$ as $k \to \infty$ for some $c_1 \in (0, \infty)$.*

Simulations show that the conjecture is true and $c_1 \approx 0.36$. The weaker form $p_k \asymp k^{-1/4}$ of Conjecture 1 was proved by Sinai [22], but only for integrated symmetric Bernoulli random walks. It also interesting to note that, by McKean [19], the unilateral small deviation probabilities of an integrated Wiener process have the same order as $T \to \infty$:

$$(43) \quad \mathbb{P}\left\{\min_{0 \leq s \leq T} \int_0^s W(u)\,du \geq -1\right\} \sim c_2 T^{-1/4}$$

for some $c_2 \in (0, \infty)$. The left-hand side of (43) is a unilateral small deviation probability since

$$(44) \quad \mathbb{P}\left\{\min_{0 \leq s \leq T} \int_0^s W(u)\,du \geq -1\right\} = \mathbb{P}\left\{\min_{0 \leq s \leq 1} \int_0^s W(u)\,du \geq -T^{-3/2}\right\}.$$

To be precise, McKean was interested in a more general problem, and some calculations are required to obtain (43) from his results. Therefore, we additionally refer to Isozaki and Watanabe [12] who state (43) explicitly.

By the results mentioned above, we also suppose that Conjecture 1 is true for other integrated centered random walks that satisfy some moment conditions.

Now we are able to prove the following result on convergence of $\frac{K_n(1)}{\sqrt{n}}$.

PROPOSITION 5. *Suppose Conjecture 1 holds true. Then in the Poisson model, we have*

$$(45) \qquad \lim_{n\to\infty} \mathbb{E}\left(\frac{K_n(1)}{\sqrt{n}}\right) = c_3, \qquad \sup_{n\geq 1} \mathbb{E}\left(\frac{K_n(1)}{\sqrt{n}}\right)^2 < \infty$$

*for some $c_3 \in (0,\infty)$; the sequence $\frac{K_n(1)}{\sqrt{n}}$ is tight and uniformly integrable; and the limit of any weakly converging subsequence of $\frac{K_n(1)}{\sqrt{n}}$ takes value zero with positive probability, but is not identically equal to zero.*

Numerical simulations show that $\frac{K_n(1)}{\sqrt{n}}$ is weakly convergent and that this convergence is quite fast. In Figure 1 we present the (empirical) distribution function of $\frac{K_n(1)}{\sqrt{n}}$ for $n = 10{,}000$. Since the simulations performed for $n = 40{,}000$ showed a hardly perceptible difference, this function seems to be a good candidate for the distribution function of the conjectured limit.

Note that if we weaken Conjecture 1 to $p_k \asymp k^{-1/4}$, then Proposition 5 still holds true with the only difference that $\mathbb{E}\frac{K_n(1)}{\sqrt{n}} \asymp 1$.

PROOF OF PROPOSITION 5. We start with the convergence of the expectation. On the one hand, by (6) and Lemma 2,

$$\mathbb{E}\left(\frac{K_n(1)}{\sqrt{n}}\right) = \frac{1}{\sqrt{n}} + \frac{e}{\sqrt{n}} \sum_{i=1}^{n-1} p_i p_{n-i},$$

and on the other hand,

$$\frac{1}{\sqrt{n}} \sum_{i=1}^{n-1} i^{-1/4}(n-i)^{-1/4} = \frac{1}{n} \sum_{i=1}^{n-1} \left(\frac{i}{n}\right)^{-1/4} \left(1 - \frac{i}{n}\right)^{-1/4} \longrightarrow \mathrm{B}(3/4, 3/4)$$

as the integral sum of Beta function. Then it follows from Conjecture 1 and standard arguments that $\mathbb{E}\frac{K_n(1)}{\sqrt{n}}$ converges to $c_3 := ec_1^2 \mathrm{B}(3/4, 3/4) > 0$.

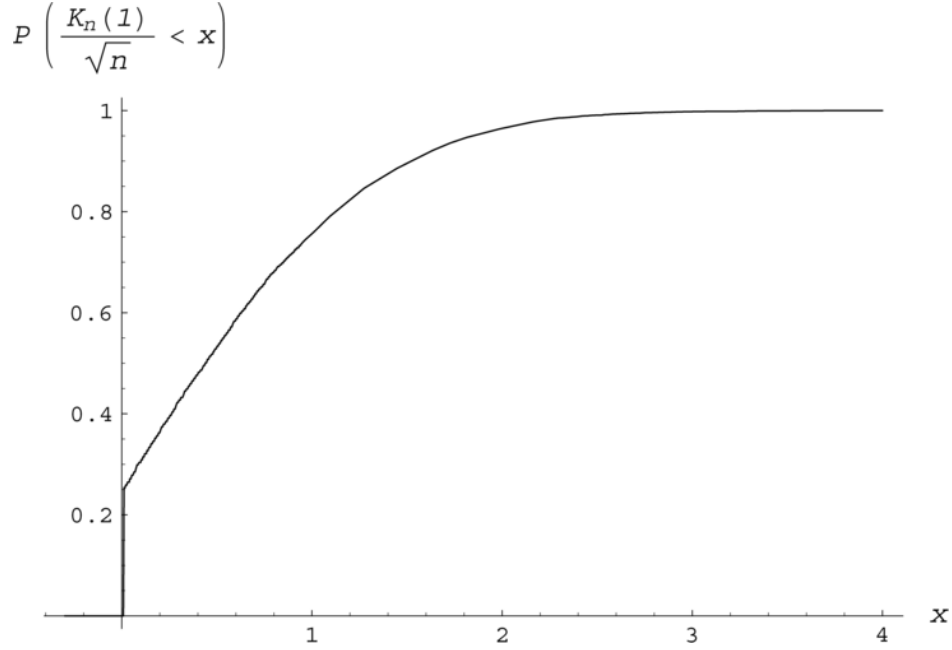

Fig. 1. *The distribution function of $\frac{K_n(1)}{\sqrt{n}}$ for $n = 10,000$.*

Now we check the uniform boundedness of $\mathbb{E}(\frac{K_n(1)}{\sqrt{n}})^2$. By (6) it is sufficient to prove that

(46) $$\sup_{n \geq 1} \frac{1}{n} \sum_{i,j=1, i \neq j}^{n-1} \mathbb{P}\{T_{i,n} > 1, T_{j,n} > 1\} < \infty.$$

Suppose $i < j$; then by using (8) and properties of $F_{k,j,l}(\cdot)$, we get

$$\mathbb{P}\{T_{i,n} > 1, T_{j,n} > 1\} = \mathbb{P}\left\{\min_{\substack{1 \leq k \leq n-i \\ 1 \leq l \leq i}} F_{k,i,l}(1) > 0, \min_{\substack{1 \leq k \leq n-j \\ 1 \leq l \leq j}} F_{k,j,l}(1) > 0\right\}$$

$$\leq \mathbb{P}\left\{\min_{\substack{1 \leq k \leq (j-i)/2 \\ 1 \leq l \leq i}} F_{k,i,l}(1) > 0, \min_{\substack{1 \leq k \leq n-j \\ 1 \leq l \leq (j-i)/2}} F_{k,j,l}(1) > 0\right\},$$

where by $(j-i)/2$ we mean $\lceil (j-i)/2 \rceil$. The minima in the last expression are independent as functions of $\{X_m\}_{m \leq (i+j)/2}$ and $\{X_m\}_{m \geq (i+j)/2+1}$, respectively; hence

$$\mathbb{P}\{T_{i,n} > 1, T_{j,n} > 1\} \leq \mathbb{P}\left\{\min_{\substack{1 \leq k \leq (j-i)/2 \\ 1 \leq l \leq i}} F_{k,i,l}(1) > 0\right\} \cdot \mathbb{P}\left\{\min_{\substack{1 \leq k \leq n-j \\ 1 \leq l \leq (j-i)/2}} F_{k,j,l}(1) > 0\right\}$$

$$= \mathbb{P}\{T_{i,i+(j-i)/2} > 1\} \cdot \mathbb{P}\{T_{(j-i)/2, n-j+(j-i)/2} > 1\}$$
$$= e^2 p_i p_{\lceil(j-i)/2\rceil}^2 p_{n-j},$$

where the first equality follows from (8) and the second follows from Lemma 2.

Recalling Conjecture 1, we get

$$\frac{1}{n} \sum_{i,j=1, i\neq j}^{n-1} \mathbb{P}\{T_{i,n} > 1, T_{j,n} > 1\} \leq \frac{1}{n} \sum_{i,j=1, i\neq j}^{n-1} e^2 p_i p_{\lceil|j-i|/2\rceil}^2 p_{n-j}$$

$$\leq \frac{c}{n} \sum_{i,j=1, i\neq j}^{n-1} i^{-1/4} \lceil|j-i|/2\rceil^{-1/2} (n-j)^{-1/4}$$

$$\leq \frac{c}{n^2} \sum_{i,j=1, i\neq j}^{n-1} \left(\frac{i}{n}\right)^{-1/4} \left|\frac{j}{n} - \frac{i}{n}\right|^{-1/2} \left(1 - \frac{j}{n}\right)^{-1/4}$$

for some $c > 0$. The last expression is an integral sum converging to

$$c \int_0^1 \int_0^1 x^{-1/4} |x-y|^{-1/2} (1-y)^{-1/4} \, dx \, dy,$$

and it is a simple exercise to check that the integral is finite. This concludes (46).

The uniform integrability of $\frac{K_n(1)}{\sqrt{n}}$ follows from the second relation from (45), see Billingsley [3], Section 5, and the tightness follows from the uniform integrability.

Finally, suppose $\frac{K_{n_i}(1)}{\sqrt{n_i}} \xrightarrow{\mathcal{D}} \xi$ for some subsequence $n_i \to \infty$ and some r.v. $\xi$. Then $\mathbb{E}\xi = c_3 > 0$ by the uniform integrability and (45), and hence $\xi$ is not identically equal to zero. But the distribution of $\xi$ has an atom at zero since by Proposition 4 and properties of weak convergence,

$$\mathbb{P}\{\xi = 0\} = \lim_{\varepsilon \searrow 0} \mathbb{P}\{\xi \leq \varepsilon\}$$

$$\geq \lim_{\varepsilon \searrow 0} \limsup_{i \to \infty} \mathbb{P}\left\{\frac{K_{n_i}(1)}{\sqrt{n_i}} \leq \varepsilon\right\}$$

$$\geq \lim_{\varepsilon \searrow 0} \lim_{i \to \infty} \mathbb{P}\{K_{n_i}(1) = 1\} > 0. \qquad \square$$

**7. Open questions.** 1. The number of clusters at the critical moment $t = 1$.

Here the main question is if Conjecture 1 holds true. Even by itself, this problem is worth studying.

But even if Conjecture 1 is true, we still do not have a proof of weak convergence of $\frac{K_n(1)}{\sqrt{n}}$, it is only known that this sequence is tight. The author strongly believes, relying on numerical simulations, that the limit exists. It would be interesting to find this conjectured limit, which should be nontrivial by Proposition 5, in an explicit form.

2. The weak convergence of $\frac{K_n(\cdot)-na(\cdot)}{\sqrt{n}}$ on the whole interval $[0,1]$.

It is very natural to ask if it is possible to strengthen Theorem 1 by proving the weak convergence of $\frac{K_n(\cdot)-na(\cdot)}{\sqrt{n}}$ in $D[0,1]$. This complicated problem returns us again to Question 1 because the weak convergence of $\frac{K_n(\cdot)-na(\cdot)}{\sqrt{n}}$ in $D[0,1]$ implies the weak convergence of $\frac{K_n(1)-na(1)}{\sqrt{n}} = \frac{K_n(1)}{\sqrt{n}}$, see Billingsley [3], Section 15. But even if $\frac{K_n(1)}{\sqrt{n}}$ converges, its weak limit $K(1)$ is not Gaussian, hence the limit process $K(\cdot)$, which is Gaussian on $[0,1)$, is no more Gaussian on $[0,1]$. Therefore, it is doubtful that Theorem 1 is true in $D[0,1]$; at least, one should provide a proof completely different from the presented one. Also, it is unclear how to define the finite-dimensional distributions of the non-Gaussian $K(\cdot)$ on $[0,1]$ because simulations show that $K(1)$ would not be independent with $K(t)$ for $t<1$.

3. The number of clusters in the warm gas.

In the presented case, initial speeds of particles are zero. This model is often called the *cold* gas according to its zero initial temperature. We introduce a new model stating that initial speeds of particles are $a_n v_1, a_n v_2, \ldots, a_n v_n$, where $v_i$ are some i.i.d. r.v.'s and $a_n$ is a sequence of normalization constants. This model, called the *warm* gas, was studied in many papers, for example, [14, 16, 20, 25].

It is of a great interest to study the behavior of $K_n(t)$ in the warm gas. In [25], the author proved that in the basic case where $a_n = 1$ for all $n$ and $\mathbb{E}v_i^2 < \infty$, we have $\frac{K_n(t)}{n} \xrightarrow{\mathbb{P}} 0$ for all $t > 0$. The question is to find a normalization of $K_n(t)$ leading to some nontrivial limit. Clearly, this normalization depends on $a_n$, but it is very possible that there is an effect of phase transition similar to the one discovered by Lifshits and Shi [16]: If $a_n$ are small enough, then the gas has a low temperature and the normalization is the same as in the cold gas. If $a_n$ are big enough, as in the basic case $a_n \equiv 1$, then the normalization and the behavior of the gas differ entirely from the case of the cold gas.

The author believes that the localization property, which is described in Section 3, could be helpful in a study of these questions.

It is also interesting to compare the behavior of $K_n(1)$ in the warm and in the cold gases; in the warm gas, the moment $t=1$ plays the same "critical" role as in the cold gas, see Lifshits and Shi [16]. The variable $K_n(1)$ was studied by Suidan [24], who considered the warm gas with $a_n \equiv 1$ and deterministic initial positions of particles (his initial positions were $\frac{1}{n}, \frac{2}{n}, \ldots, \frac{n}{n}$). For this case, Suidan found the distribution of $K_n(1)$ and showed that $\mathbb{E}K_n(1) \sim \log n$. Recall that in the presented case, $\mathbb{E}K_n(1) \sim c_3 \sqrt{n}$.

CLUSTERING IN A STOCHASTIC MODEL OF ONE-DIMENSIONAL GAS 33## 4. The number of clusters in ballistic systems of sticky particles.

A sticky particles model is called *ballistic* if it evolves according to the laws introduced in Section 1, but in the absence of gravitation. Such models are, in some sense, more natural than gravitational ones because the basic assumption that gravitation does not depend on distance is sometimes confusing. However, an unpublished paper of Lifshits and Kuoza shows that certain gravitational and ballistic models are tightly connected.

It seems interesting to study the number of clusters in the ballistic model. The author does not know any results in this field.

**Acknowledgments.** I am grateful to my adviser Mikhail A. Lifshits for drawing my attention into the subject and for his guidance. I also thank the anonymous referees for carefully reading this paper and useful comments.## REFERENCES

[1] BAUM, L. E. and KATZ, M. (1965). Convergence rates in the law of large numbers. *Trans. Amer. Math. Soc.* **120** 108–123. MR0198524
[2] BERTOIN, J. (2002). Self-attracting Poisson clouds in an expanding universe. *Comm. Math. Phys.* **232** 59–81. MR1942857
[3] BILLINGSLEY, P. (1968). *Convergence of Probability Measures*. Wiley, New York. MR0233396
[4] BRENIER, Y. and GRENIER, E. (1998). Sticky particles and scalar conservation laws. *SIAM J. Numer. Anal.* **35** 2317–2328. MR1655848
[5] CHERTOCK, A., KURGANOV, A. and RYKOV, YU. (2007). A new sticky particle method for pressureless gas dynamics. *SIAM J. Numer. Anal.* **45** 2408–2441. MR2361896
[6] E, W., RYKOV, YU. G. and SINAI, YA. G. (1996). Generalized variational principles, global weak solutions and behavior with random initial data for systems of conservation laws arising in adhesion particle dynamics. *Comm. Math. Phys.* **177** 349–380. MR1384139
[7] ESARY, J. D., PROSCHAN, F. and WALKUP, D. W. (1967). Association of random variables, with applications. *Ann. Math. Stat.* **38** 1466–1474. MR0217826
[8] GIRAUD, C. (2001). Clustering in a self-gravitating one-dimensional gas at zero temperature. *J. Statist. Phys.* **105** 585–604. MR1871658
[9] GIRAUD, C. (2005). Gravitational clustering and additive coalescence. *Stochastic Process. Appl.* **115** 1302–1322. MR2152376
[10] GURBATOV, S. N., MALAKHOV, A. N. and SAICHEV, A. I. (1991). *Nonlinear Random Waves and Turbulence in Nondispersive Media*: *Waves, Rays, Particles*. Manchester Univ. Press. MR1255826
[11] GURBATOV, S. N., SAICHEV, A. I. and SHANDARIN, S. F. (1989). The large-scale structure of the universe in the frame of the model equation of nonlinear diffusion. *Mon. Not. R. Astr. Soc.* **236** 385–402.
[12] ISOZAKI, Y. and WATANABE, S. (1994). An asymptotic formula for the Kolmogorov diffusion and a refinement of Sinai's estimates for the integral of Brownian motion. *Proc. Japan Acad. Ser. A Math. Sci.* **70** 271–276. MR1313176
[13] KARLIN, S. (1968). *A First Course in Stochastic Processes*. Academic Press, New York. MR0208657

DEPARTMENT OF PROBABILITY THEORY
AND MATHEMATICAL STATISTICS
FACULTY OF MATHEMATICS AND MECHANICS
ST. PETERSBURG STATE UNIVERSITY
BIBLIOTECHNAYA PL. 2
STARY PETERHOF 198504
RUSSIA
E-MAIL: vlad.vysotsky@gmail.com